\documentclass{article}
\pdfoutput=1

\newcommand{\sect}[1]{\section{#1}\setcounter{equation}{0}}

\newtheorem{thh}{Theorem}[section]
\newtheorem{cor}[thh]{Corollary}
\newtheorem{lem}[thh]{Lemma}
\newtheorem{prop}[thh]{Proposition}

\def\sep{\;\vrule\;}
\def\proof#1. {\par
                      \ifdim\lastskip<15pt
                      \removelastskip\penalty-200
                      \vskip15pt plus3pt minus3pt
                      \fi
                       {\def\a{#1}
                       \ifx\a\empty
                       {\noindent\bf Proof.}
                       \else
                       {\noindent\bf Proof of #1.}
                       \fi}\enspace}
\def\restr#1{\,\vrule\,\lower1.75ex\hbox{$#1$}}

\def\endproof{\hfill\hspace{-6pt}\rule[-14pt]{6pt}{6pt}
\vskip22pt plus3pt minus 3pt}

\def\be{\begin{equation}}
\def\ee{\end{equation}}
\def\bea{\begin{eqnarray}}
\def\eea{\end{eqnarray}}
\def\bean{\begin{eqnarray*}}
\def\eean{\end{eqnarray*}}

\def\a{\alpha}
\def\b{\beta}
\def\d{\delta}

\def\e{\varepsilon}
\def\f{\varphi}
\def\F{\Phi}
\def\g{\gamma}
\def\G{\Gamma}
\def\i{\infty}

\def\O{\Omega}

\def\s{\sigma}

\def\t{\theta}

\def\z{\zeta}

\def\c{{\rm cap}}

\def\setm{\setminus}

\def\ov{\overline}

\def\c{{\rm cap}}

\font\tenopen = cmbx10
\font\sevenopen = cmbx7
\font\fiveopen = cmbx5
\newfam\openfam
\def\open{\fam\openfam\tenopen}
\textfont\openfam = \tenopen
\scriptfont\openfam = \sevenopen
\scriptscriptfont\openfam = \fiveopen

\def\C{\overline{\open C}}

\newcommand*{\unitcircle}{\mathbb{T}}

\newcommand*{\innerdisk}{D_1}
\newcommand*{\outercirc}{\mathbb{T}(r_1)}
\newcommand*{\outerdisk}{D_2}

\newcommand*{\innercurve}{\Phi_2(\partial \outerdisk)}
\newcommand*{\outercurve}{\Phi_1(\partial \innerdisk)}

\usepackage{graphicx}
\usepackage{amsmath}
\usepackage{amsfonts}

\title{Bernstein- and Markov-type inequalities for rational functions}
\author{Sergei Kalmykov,  B\'ela Nagy and Vilmos Totik}
\date{}
\begin{document}
\maketitle
\begin{abstract} Asymptotically sharp Bernstein- and Markov-type
inequalities are established for rational functions on $C^2$ smooth
Jordan curves and arcs.
The results are formulated in terms of the normal derivatives of certain Green's functions with poles
at the poles of the rational functions in question. As a special case
(when all the poles
are at infinity) the corresponding results for polynomials are recaptured.
\end{abstract}

\tableofcontents

\sect{Introduction}
Inequalities for polynomials  have a rich
history and numerous applications in different branches of mathematics,
in particular in approximation theory (see, for example, the books
\cite{BE}, \cite{DevoreLorentz} and \cite{Milovanovic}, as well as  the
extensive references there).
The two most classical results are the
Bernstein inequality \cite{Bernstein}
\be|P_n'(x)|\le \frac{n}{\sqrt{1-x^2}}\|P_n\|_{[-1,1]},\qquad x\in (-1,1),\label{bb1}\ee
and the Markov inequality \cite{Markov}
\be \|P_n'\|_{[-1,1]}\le n^2\|P_n\|_{[-1,1]}\label{m1}\ee
for estimating the derivative of polynomials $P_n$ of degree at most
$n$ in terms of
the supremum norm $\|P_n\|_{[-1,1]}$ of the polynomials. In (\ref{bb1}) the
order of the right hand side  is $n$, and the estimate can be used at inner points
of $[-1,1]$. In (\ref{m1}) the growth of the right-hand side
is $n^2$, which is much larger, but  (\ref{m1}) can also be used close to the endpoints
$\pm 1$, and it gives a global estimate. We shall use the terminology ``Bernstein-type
inequality" for estimating the derivative away from endpoints with a
factor of order $n$, and ``Markov-type inequality" for a global estimate
on the derivative with a factor of order $n^2$.

The Bernstein and Markov inequalities have been generalized and improved in several directions
over the last century, see the extensive books \cite{BE}
and \cite{Milovanovic}.
See also \cite{Dubinin} and the references there for various improvements.
For rational functions sharp Bernstein-type inequalities have been given
for circles \cite{BorweinErdelyipaper} and for compact subsets of the real line and
 circles, see \cite{BorweinErdelyipaper}, \cite{DubininKalmykov}, \cite{Lukashov2004}.
 We are
 unaware of a corresponding Markov-type estimate.
General (but not sharp) estimates on the derivative of rational functions can also be found in
\cite{Ruszak1} and \cite{Ruszak2}.

The aim of this paper is to give the sharp form
of the Bernstein and Markov inequalities for rational functions on
smooth Jordan curves and arcs. We shall be primarily interested in
the asymptotically best possible estimates and in the structure of the
constants on the right hand side. As we shall see,
from this point of view there is a huge difference
between Jordan curves and Jordan arcs.
All the results are formulated in terms of the normal derivatives of certain Green's functions with poles
at the poles of the rational functions in question. When all the poles
are at infinity we recapture the corresponding results for polynomials
that have been proven in the last decade.

We shall use basic notions of potential theory, for the necessary background  we refer to
the books
\cite{Gardiner}, \cite{Ransford}, \cite{SaffTotik} or \cite{Tsuji}.

\sect{Results}\label{sectresults}

We shall work with Jordan curves and Jordan arcs on the plane. Recall that  a Jordan curve
is a homeomorphic image of a circle, while a Jordan arc is a homeomorphic image
of a segment. We say that the Jordan arc $\G$ is $C^2$ smooth if it
has a parametrization $\g(t)$, $t\in [-1,1]$, which is twice continuously
differentiable and $\g'(t)\not=0$ for $t\in [-1,1]$. Similarly we speak of $C^2$ smoothness
of a Jordan curve, the only difference is that for a Jordan curve the parameter domain
is the unit circle.

If $\G$ is a Jordan curve, then we think it counterclockwise oriented. $\C\setm \G$ has
two connected components, we denote the bounded component by $G_-$ and the unbounded
one by $G_+$. At a point $z\in \G$ we denote the two normals to $\G$ by ${\bf n}_\pm=
{\bf n}_\pm(z)$ with the agreement that ${\bf n}_-$ points towards $G_-$. So, as we
move on $\G$ according to its orientation, ${\bf n}_-$ is the left
and ${\bf n}_+$ is the right normal. In a similar fashion, if $\G$ is a Jordan
arc then we take an orientation of $\G$ and let
${\bf n}_-$ resp.  ${\bf n}_+$ denote the left resp. right normal to
$\G$ with respect to this orientation.

Let $R$ be a rational function. We say it has total degree $n$ if the sum of
the order of its poles  (including the possible pole at $\i$) is $n$. We shall often
use summations $\sum_a$ where $a$ runs through the poles of $R$, and let us agree
that in such sums a pole $a$ appears as many times as its order.

In this paper
we determine the asymptotically sharp
analogues of the Bernstein and Markov inequalities on Jordan
curves and arcs $\G$ for rational functions. Note however, that even in the simplest
case $\G=[-1,1]$ there is no Bernstein- or Markov-type inequality just in terms of
the degree of the rational function. Indeed, if $M>0$, then $R_2(z)=1/(1+Mz^2)$ is at most
1 in absolute value on $[-1,1]$, but $|R_2'(1/\sqrt M)|= \sqrt M/2$, which can
be arbitrary large if $M$ is large. Therefore, to get
Bernstein-Markov-type inequalities in the classical sense we
should limit the poles of $R$ to lie far from $\G$. In this paper
we assume that the poles of the rational functions lie in
a closed set $Z\subset \C\setm \G$ which we fix in advance.
If $Z=\{\i\}$, then $R$ has to be a polynomial.

In what follows $\|f\|_\G=\sup_{z\in \G}|f(z)|$ denotes the supremum
norm on $\G$, and $g_G(z,a)$ the Green's function of a domain $G$
with pole at $a\in G$.

Our first result is a Bernstein-type inequality on Jordan curves.

\begin{thh}\label{thcurve} Let $\G$ be a $C^2$ smooth Jordan curve on the plane, and let
$R_n$ be a rational function of total degree $n$ such that
its poles lie in the fixed closed set $Z\subset \C\setm \G$. If $z_0\in \G$, then
\be |R_n'(z_0)|\le (1+o(1)) \|R_n\|_{\G}\max\left(\sum_{a\in Z\cap G_+}\frac{\partial g_{G_+}(z_0,a)}{\partial {\bf n}_+},
\sum_{a\in Z\cap G_-}\frac{\partial g_{G_-}(z_0,a)}{\partial {\bf n}_-}
\right),\label{maincurve}\ee
where the summations are for the poles of $R_n$ and where
$o(1)$ denotes a quantity that tends to 0 uniformly in $R_n$ as $n\to\i$.
Furthermore, this estimate holds uniformly in $z_0\in \G$.
\end{thh}
The normal derivative $\partial g_{G_\pm}(z_0,a)/\partial {\bf n}_\pm$ is $2\pi$-times the
density of the harmonic measure of $a$ in the domain $G_\pm$, where
the density is taken with respect to the arc measure on $\G$. Thus, the right
hand side in (\ref{maincurve}) is easy to formulate in terms
of harmonic measures, as well.

\begin{cor}\label{corthcurve} If $\G$ is as in Theorem \ref{thcurve} and $P_n$ is a polynomial
of degree at most $n$, then for $z_0\in \G$ we have
\be |P_n'(z_0)|\le (1+o(1)) n\|P_n\|_{\G}\frac{\partial g_{G_+}(z_0,\i)}{\partial {\bf n}_+}.\label{maincurvecor}\ee
\end{cor}
This is Theorem 1.3 in the paper \cite{NagyTotik1}. The estimate
(\ref{maincurvecor}) is asymptotically the
best possible (see below), and on the right
$\partial g_{G_+}(z_0,\i)/\partial {\bf n}_+$ is $2\pi$-times of the
density of the equilibrium measure of $\G$ with respect to the
arc measure on $\G$. Therefore,  the corollary shows an explicit relation in
between the Bernstein factor at a given point and the harmonic density at the
same point.

If $R_n$ has order $n+o(n)$ and we take the sum on the right of
(\ref{maincurve}) only on some of its $n$ poles, then (\ref{maincurve}) still holds
(i.e. $o(n)$ poles do not have to be accounted for). Now in this sense
Theorem \ref{thcurve} is sharp.
\begin{thh}\label{thcurvesh} Let $\G$ be as in Theorem \ref{thcurve} and let
$Z\subset \C\setm \G$ be a non-empty closed set. If $\{a_{1,n},\ldots,a_{n,n}\}$, $n=1,2,\ldots$
is an array of points from $Z$ and $z_0\in \G$ is a point
on $\G$, then there are non-zero rational functions
$R_n$ of degree $n+o(n)$ such that $a_{1,n},\ldots,a_{n,n}$ are among the poles of $R_n$ and
\be |R_n'(z_0)|\ge (1-o(1)) \|R_n\|_{\G}\max\left(\sum_{a_{j,n}\in G_+}\frac{\partial g_{G_+}(z_0,a_{j,n})}{\partial {\bf n}_+},
\sum_{a_{j,n}\in G_-}\frac{\partial g_{G_-}(z_0,a_{j,n})}{\partial {\bf n}_-}
\right).\label{maincurvesh}\ee
\end{thh}
In this theorem if a point $a\in Z$ appears $k$ times in $\{a_{1,n},\ldots,a_{n,n}\}$,
then the understanding is that at $a$ the rational function  $R_n$  has a pole of order $k$.

Next, we consider the Bernstein-type inequality for rational functions on a Jordan arc.

\begin{thh}\label{tharc} Let $\G$ be a $C^2$ smooth Jordan arc on the plane, and let
$R_n$ be a rational function of total degree $n$
such that
its poles lie in the fixed closed set $Z\subset \C\setm \G$. If $z_0\in \G$ is different
from the endpoints of $\G$, then
\be |R_n'(z_0)|\le (1+o(1)) \|R_n\|_{\G}
\max\left(\sum_{a\in Z}\frac{\partial g_{\C\setm \G}(z_0,a)}{\partial {\bf n}_+},
\sum_{a\in Z}\frac{\partial g_{\C\setm \G}(z_0,a)}{\partial {\bf n}_-}\right),\label{mainarc}\ee
where the summations are for the poles of $R_n$ and where
$o(1)$ denotes a quantity that tends to 0 uniformly in $R_n$ as $n\to\i$.
Furthermore, (\ref{mainarc}) holds uniformly in $z_0\in J$ for any closed subarc $J$ of $\G$ that does
not contain either of the endpoints of $\G$.
\end{thh}

\begin{cor}\label{cortharc} If $\G$ is as in Theorem \ref{tharc} and $P_n$ is a polynomial
of degree at most $n$, then for $z_0\in \G$, which is different from the endpoints
of $\G$, we have
\be |P_n'(z_0)|\le (1+o(1)) n\|P_n\|_{\G}\max\left(\frac{\partial g_{\C\setm \G}(z_0,\i)}{\partial {\bf n}_+},
\frac{\partial g_{\C\setm \G}(z_0,\i)}{\partial {\bf n}_-}\right).\label{mainarcor}\ee
\end{cor}
This was proven in \cite{NK} for analytic $\G$ and in \cite{TotikMarkov} for $C^2$ smooth $\G$.
More generally, if $a_1,\ldots,a_m$ are finitely many fixed points outside $\G$ and
\be R_n(z)=P_{n_0,0}(z)+\sum_{i=1}^mP_{n_i,i}\left(\frac{1}{z-a_i}\right)\label{spec}\ee
where $P_{n_i,i}$ are polynomials of degree at most $n_i$, then, as $n_0+\cdots +n_m\to\i$,
\be |R_n'(z_0)|\le (1+o(1)) \|R_n\|_{\G}\max\left(\sum_{i=0}^m
n_i\frac{\partial g_{\C\setm \G}(z_0,a_i)}{\partial {\bf n}_+},
\sum_{i=0}^mn_i\frac{\partial g_{\C\setm \G}(z_0,a_i)}{\partial {\bf n}_-}\right),\label{mainarc009}\ee
where $a_0=\i$.

Theorem \ref{tharc} is sharp again regarding the Bernstein factor on the right.
\begin{thh}\label{tharcsh} Let $\G$ be as in Theorem \ref{tharc} and let
$Z\subset \C\setm \G$ be a non-empty  closed set. If $\{a_{1,n},\ldots,a_{n,n}\}$, $n=1,2,\ldots$
is an arbitrary array of points from $Z$ and $z_0\in \G$ is any point
on $\G$ different from the endpoints of $\G$, then there are non-zero rational functions
$R_n$ of degree $n+o(n)$ such that $a_{1,n},\ldots,a_{n,n}$ are among the poles of $R_n$ and
\be |R_n'(z_0)|\ge (1-o(1)) \|R_n\|_{\G}
\max\left(\sum_{a\in Z}\frac{\partial g_{\C\setm \G}(z_0,a)}{\partial {\bf n}_+},
\sum_{a\in Z}\frac{\partial g_{\C\setm \G}(z_0,a)}{\partial {\bf n}_-}\right).\label{mainarcsh}\ee
\end{thh}

Now we consider the Markov-type inequality
on a $C^2$ Jordan arc $\G$ for rational functions of the form (\ref{spec}).
Let $A,B$ be the two endpoints of $\G$.
We need the quantity
\be \O_a(A)=\lim_{z\to A, \ z\in \G}\sqrt{|z-A|}\frac{\partial g_{\C\setm \G}(z,a)}{\partial {\bf n}_\pm(z)}.
\label{Odef}\ee
It will turn out that this limit exists and it is the same if we use in it the left
or the right  normal derivative (i.e. it is indifferent if we use
${\bf n}_+$  or ${\bf n}_-$ in the definition). We define $\O_a(B)$ similarly.
With these we have

\begin{thh} \label{thMarkov} Let $\G$ be a $C^2$ smooth Jordan arc on the plane, and let
$R_n$ be a rational function of total degree $n$ of the form (\ref{spec}) with fixed
$a_0,a_1,\ldots,a_m$. Then
\be \|R_n'\|_\G\le (1+o(1)) \|R_n\|_{\G}
2\max\left(\sum_{i=0}^m n_i\O_{a_i}(A),
\sum_{i=0}^m n_i\O_{a_i}(B)\right)
^2,\label{mainMarkov}\ee
where $o(1)$ tends to 0 uniformly in $R_n$ as $n\to\i$.
\end{thh}

Theorem \ref{thMarkov} is again the best possible, but we shall not state
that, for we will have
a more general result in Theorem \ref{thMarkovk}.

Actually, there is a separate Markov-type inequality around both endpoints  $A$ and $B$.
Indeed, let $U$ be  a closed neighborhood of $A$ that does not contain $B$. Then
\be \|R_n'\|_{\G\cap U}\le (1+o(1)) \|R_n\|_{\G}
2\left(\sum_{i=0}^m n_i\O_{a_i}(A)\right)^2,\label{mainMarkovalt}\ee
and this is sharp. Now (\ref{mainMarkov})
is clearly a consequence of this  and its analogue for the endpoint
$B$. Note that the discussion below will show that  the right-hand side in
(\ref{mainMarkov}) is of size $\sim n^2$, while on any
closed Jordan subarc of $\G$ that does not contain $A$ or $B$
 the derivative $R_n'$ is $O(n)$.

Let us also mention that in these theorems in general the $o(1)$ term in the $1+o(1)$ factors on
the right cannot be omitted. Indeed, consider for example, Corollary \ref{corthcurve}.
It is easy to construct a $C^2$ Jordan curve for which the normal derivative
on the right of (\ref{maincurvecor}) is small, so with   $P_1(z)=z$ the inequality
in (\ref{maincurvecor}) fails if we write $0$ instead of $o(1)$.
\bigskip

It is also interesting to consider higher derivatives, though we can do
a complete analysis only for rational functions of the form
(\ref{spec}). For them the inequalities
(\ref{maincurve}) and (\ref{mainarc}) can simply be
iterated. For example,  if $\G$ is a Jordan  arc, then under the assumptions of Theorem \ref{tharc} we have
for any fixed $k=1,2,\ldots$
\be |R_n^{(k)}(z_0)|\le (1+o(1)) \|R_n\|_{\G}
\max\left(\sum_{i=0}^m
n_i\frac{\partial g_{\C\setm \G}(z_0,a_i)}{\partial {\bf n}_+},
\sum_{i=0}^mn_i\frac{\partial g_{\C\setm \G}(z_0,a_i)}{\partial {\bf n}_-}\right)^k\label{mainarck}\ee
uniformly in $z_0\in J$ where $J$ is any closed subarc of $\G$ that does not contain
the endpoints of $\G$.
It can also be proven that this inequality is sharp for every $k$ and every $z_0\in \G$
in the sense given in Theorems \ref{thcurvesh} and \ref{tharcsh}.

The situation is different for the Markov inequality (\ref{mainMarkov}), because if we iterate
it, then we do not obtain the sharp inequality for the norm of the $k$-th derivative
(just like the iteration of the classical A. A. Markov inequality does
not give the sharp V. A. Markov inequality for higher derivatives of polynomials).
Indeed, the sharp form is given in the following theorem.

\begin{thh} \label{thMarkovk} Let $\G$ be a $C^2$ smooth Jordan arc on the plane, and let
$R_n$ be a rational function of total degree $n$ of the form (\ref{spec})
with fixed $a_0,a_1,\ldots,a_m$. Then for any fixed $k=1,2,\ldots$ we have
\be \|R_n^{(k)}\|_\G\le (1+o(1)) \|R_n\|_{\G}
\frac{2^k}{(2k-1)!!}
\max\left(\sum_{i=0}^m n_i\O_{a_i}(A),
\sum_{i=0}^m n_i\O_{a_i}(B)\right)
^{2k},
\label{mainMarkovk}\ee
where $o(1)$ tends to 0 uniformly in $R_n$ as $n\to\i$.
Furthermore, this is sharp, for one cannot write
 a constant smaller than 1 instead of $1+o(1)$ on the right.
\end{thh}

Recall that $(2k-1)!!=1\cdot 3\cdot\cdots \cdot(2k-3)\cdot (2k-1)$.

As before, this theorem will follow if we prove for any closed neighborhood $U$ of the endpoint $A$ that does not
contain the other endpoint $B$ the estimate
\be \|R_n^{(k)}\|_{\G\cap U}\le (1+o(1)) \|R_n\|_{\G}
\frac{2^k}{(2k-1)!!}\left(\sum_{i=0}^m n_i\O_{a_i}(A)\right)^{2k}.\label{mainMarkovaltk}\ee

\begin{cor}\label{corthMarkovk} If $\G$ is as in Theorem \ref{thMarkovk} and
$P_n$ is a polynomial of degree at most $n$, then
\be \|P_n^{(k)}\|_\G\le (1+o(1)) \|P_n\|_{\G}
\frac{2^k}{(2k-1)!!}n^{2k}\max\left(\O_\i(A),\O_\i(B)\right)^{2k}.\label{mainMarkovkcor}\ee
\end{cor}

This was proven in \cite[Theorem 2]{TotikMarkov}.
\bigskip

The outline of the paper is as follows.
\begin{itemize}
\item After some preparations first we verify Theorem \ref{thcurve} (Bernstein-type inequality) for {\it analytic curves}
via conformal maps onto the unit disk and using on the unit
disk a result of Borwein and Erd\'elyi. This part uses in an essential
way a decomposition theorem  for meromorphic functions.

\item Next, Theorem \ref{tharc} is verified for {\it analytic
arcs} from the analytic case of Theorem \ref{thcurve} for Jordan curves
 via the Joukowskii mapping.
\item For $C^2$ arcs Theorem \ref{tharc} follows from its version for analytic
arcs by an appropriate approximation.
\item For $C^2$ curves Theorem \ref{thcurve} will be deduced
from Theorem \ref{tharc} by introducing a gap (omitting a small part) on the given
Jordan curve to get a Jordan arc, and then by closing up that gap.
\item The Markov-type inequality Theorem \ref{thMarkovk} is deduced from
the Bernstein-type inequality on arcs (Theorem \ref{tharc}, more precisely
from its higher derivative variant (\ref{mainarck})) by a symmetrization
technique during which the given endpoint where we consider the Markov-type inequality
is mapped into an inner point of a different Jordan arc.
\item Finally, in Section \ref{sectsharp}  we prove the sharpness of the theorems
using conformal maps and sharp forms of Hilbert's lemniscate theorem.
\end{itemize}

\sect{Preliminaries}\label{sectprelim}
In this section we collect some tools that are
used at various places in the proofs.

\subsection{A ``rough'' Bernstein-type inequality}

We  need the following ``rough'' Bernstein-type inequality on
Jordan curves.
\begin{prop}
\label{prop:rough_berns}
Let $\Gamma$ be a $C^{2}$ smooth Jordan curve and $Z\subset\C\setminus\Gamma$
a closed set.
Then there exists $C >0$ such that for any
rational function $R_n$ with poles in $Z$ and
of degree $n$, we have
\[
\|R_n'\|_\G
\le
C n\left\Vert R_n\right\Vert_{\Gamma}.
\]
\end{prop}
\proof. Recall that $G_-$ denotes the inner, while $G_+$ denotes the outer
domain to $\G$. We shall need the following Bernstein-Walsh-type estimate:
\be
\left|R_n\left(z\right)\right|
\le\left\Vert R_n\right\Vert_{\Gamma}
\exp\left(\sum_{a\in Z\cap G_\pm}g_{G_\pm}\left(z,a\right)\right)
\label{BW}\ee
where the summation  is taken for $a\in Z\cap G_+$ if $z\in G_+$ (and
then $g_{G_+}$ is used) and
for $a\in Z\cap G_-$ if $z\in G_-$. Indeed, suppose, for example,
that $z\in G_-$. The function
\[\log |R_n(z)|-\left(\sum_{a\in Z\cap G_- }g_{G_-}\left(z,a\right)\right)\]
is subharmonic in $G_-$ and has boundary values $\le \log\|R_n\|_\G$ on $\G$,
so (\ref{BW}) follows from the maximum principle for subharmonic functions.

Let $z_0\in \G$ be arbitrary.
It follows from Proposition \ref{propgreencont} below that there is a $\d>0$ such that
 for dist$(z,\G)<\d$ we have for all $a\in Z$ the bound
\[g_{G_{\pm}}\left(z,a\right)
\le
C_1\mathrm{dist}\left(z,\Gamma\right)\le C_1\left|z-z_{0}\right|\]
with some constant $C_1$.

Let $C_{1/n}(z_0):=\left\{ z\sep \left|z-z_{0}\right|=1/n \right\} $
be the circle about $z_0$ of radius $1/n$ (assuming $n>2/\delta$). For $z\in C_{1/n}(z_0)$
the sum on the right of (\ref{BW}) can be bounded as
\[
\sum_{a\in Z\cap G_+}g_{G_{+}}\left(z,a\right)
\le
nC_1\left|z-z_{0}\right|\le C_1\]
if $z\in G_+$, and a similar estimate holds if  $z\in G_-$.
Therefore, $|R_n(z)|\le e^{C_1}\|R_n\|_\G$.

Now we apply Cauchy's integral  formula
\[\left|R_n'\left(z_{0}\right)\right|
=
\left|\frac{1}{2\pi i}
\int_{C_{1/n}(z_0)}\frac{R_n\left(z\right)}{\left(z-z_{0}\right)^{2}}dz\right|
 \le
\frac{1}{2\pi}\frac{2\pi}{n}\frac{\left\Vert
R_n\right\Vert_{\Gamma}e^{C_1}}
{n^{-2}}
\\=
\left\Vert R_n\right\Vert_{\Gamma}ne^{C_1},\]
which proves the proposition.
\endproof

\subsection{Conformal mappings onto the inner and outer domains}\label{sectconf}

Denote by $\mathbb{D}=\left\{ v\sep \ \left|v\right|<1\right\} $
the unit disk and by $\mathbb{D}_+=\left\{ v\sep \ \left|v\right|>1\right\} \cup\left\{
\infty\right\} $ its exterior.

By the Kellogg-Warschawski
theorem (see e.g. \cite[Theorem 3.6]{Pommerenke}), if $\G$
is $C^{2}$ smooth, then Riemann mappings from
$\mathbb{D},\mathbb{D}_+$
onto $G_{-}, G_{+}$, respectively, as well as their derivatives can be extended
continuously to the boundary $\G$.
Under analyticity assumption, the corresponding Riemann mappings
have extensions to larger domains. In fact, the following proposition holds
(see e.g.  Proposition 7 in \cite{NK} with slightly
different notation).

\begin{prop}
\label{prop:green_deriv_g_one_g_two}Assume that $\G$ is analytic, and let $z_{0}\in \G$
be fixed. Then there exist two Riemann mappings
$\Phi_{1}:\mathbb{D}\rightarrow G_{-}$,
$\Phi_{2}:\mathbb{D}_+\rightarrow G_{+}$ such that
$\Phi_{j}\left(1\right)=z_{0}$
and $\left|\Phi_{j}'\left(1\right)\right|=1$, $j=1,2$.
Furthermore,  there exist
$0\le r_{2}<1<r_{1}\le\infty$ such that $\Phi_{1}$ extends to
a conformal map of $D_{1}:=\left\{ v\sep\ \left|v\right|<r_{1}\right\}$
and  $\Phi_{2}$ extends to a conformal map of
$D_{2}:=\left\{ v\sep\ \left|v\right|>r_{2}\right\} \cup\left\{
\infty\right\} $.
\end{prop}

Since the argument of $\F_j'(1)$ gives the angle of the tangent line
to $\G$ at $z_0$, the arguments of $\F_1'(1)$ and of $\F_2'(1)$ must
be the same, which combined with $\left|\Phi_1'(1)\right|=
\left|\Phi_2'(1)\right|=1$ yields $\Phi_1'(1)=\Phi_2'(1)$.
Therefore,
\begin{equation}
\Phi_1(1)=\Phi_2(1)=z_0,\
\Phi_1'(1)=\Phi_2'(1),\
\left|\Phi_1'(1)\right|=
\left|\Phi_2'(1)\right|=1.
\label{eq:phione_phitwo}
\end{equation}

\begin{figure}
\begin{center}
\includegraphics[keepaspectratio,width=\textwidth]{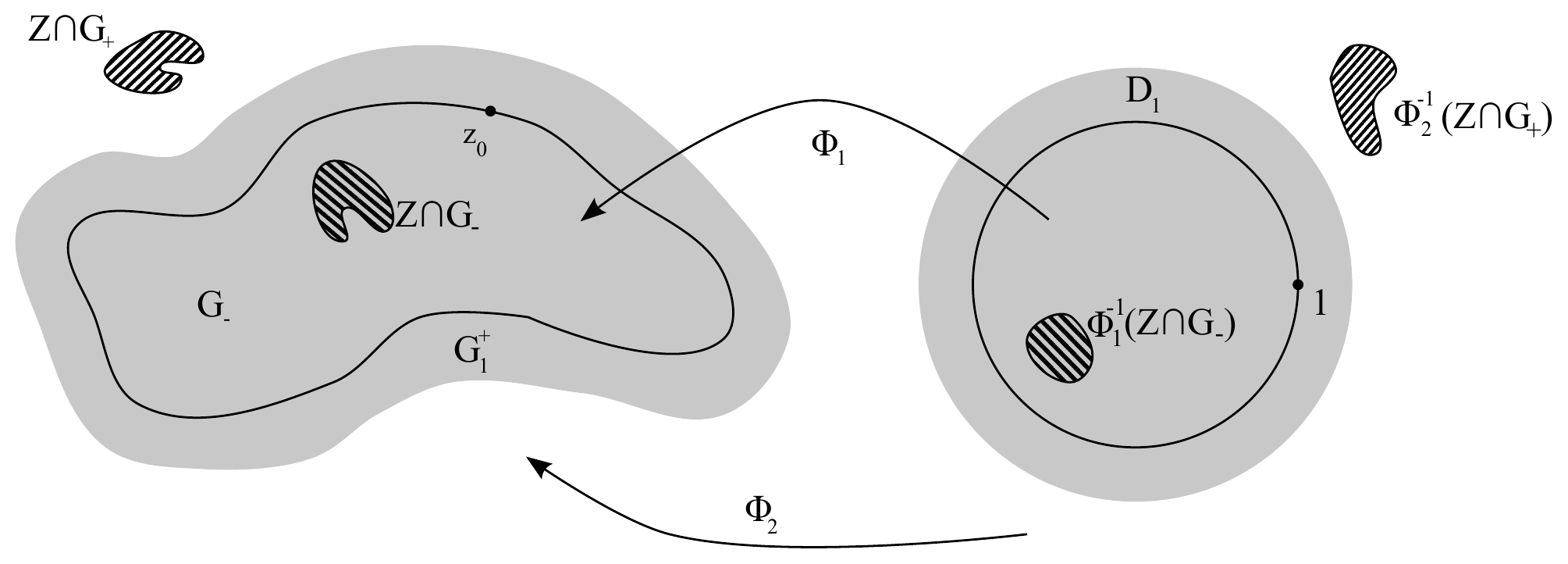}
\end{center}
\caption{The two conformal mappings $\Phi_1$, $\Phi_2$, the domain
$D_1$ and the possible location of poles}
\label{fig:analext}
\end{figure}

From now on, for a given $z_0\in \G$
we fix these two conformal maps. These mappings and
the corresponding domains are depicted on Figure \ref{fig:analext}.
We may assume that $D_1$ and $\F_2^{-1}(Z)\cap G_+$, as well as  $D_2$ and $\F_1^{-1}(Z)\cap G_-$
are of positive distance from
one another (by slightly decreasing $r_1$ and increasing $r_2$, if necessary).

\begin{prop}
\label{prop:green_deriv_unitdisk}
The following hold for arbitrary $a\in G_-$, $b \in G_+$
with $a':=\Phi_1^{-1}(a)$, $b':=\Phi_2^{-1}(b)$
\begin{gather*}
\frac{\partial g_{G_{-}}\left(z_{0},a\right)}{\partial \mathbf{n}_-}
=
\frac{\partial g_{\mathbb{D}}\left(1,a'\right)}{\partial \mathbf{n}_-}
=
\frac{1-\left|a'\right|^{2}}{\left|1-a'\right|^{2}},
\\
\frac{\partial g_{G_{+}}\left(z_{0},b\right)}{\partial \mathbf{n}_+}
=
\frac{\partial g_{\mathbb{D}_+}\left(1,b'\right)}{\partial \mathbf{n}_+}
=
\frac{\left|b'\right|^{2}-1}{\left|1-b'\right|^{2}},\qquad \mbox{if $b'\ne\infty$},
\end{gather*}
and if $b'=\infty$, then
\[
\frac{\partial g_{G_{+}}\left(z_{0},b\right)}{\partial \mathbf{n}_+}
=
\frac{\partial g_{\mathbb{D}_+}\left(1,\infty\right)}{\partial \mathbf{n}_+}
=1.
\]
\end{prop}

This proposition is a slight generalization of Proposition 8 from
\cite{NK} with the same proof.

\subsection{The Borwein-Erd\'elyi inequality}

The following inequality will be central in establishing Theorem \ref{thcurve}
in  the analytic
case, it serves as a model.
For a proof we refer to  \cite{BorweinErdelyipaper} (see also \cite[Theorem 7.1.7]{BE}).

Let $\mathbb{T}$ denote the unit circle.

\begin{prop}
[Borwein-Erd\'elyi]\label{bepr} Let $a_{1},\dots,a_{m}\in\mathbf{C}\setminus\mathbb{T}$,
\[
B_{m}^{+}\left(v\right):=\sum_{ |a_{j}|>1}\frac{|a_{j}|^{2}-1}{|a_{j}-v|^{2}},
\qquad
B_{m}^{-}\left(v\right):=\sum_{ |a_{j}|<1}\frac{1-|a_{j}|^{2}}{|a_{j}-v|^{2}},
\]
and $B_{m}\left(v\right):=\max\left(B_{m}^{+}\left(v\right),\,
B_{m}^{-}\left(v\right)\right)$.
If $P$ is a polynomial with $\deg(P)\le m$ and
$R_m\left(v\right)=P\left(v\right)/\prod_{j=1}^{m}\left(v-a_{j}\right)$
is a rational function, then
\[
\left|R_m'(v)\right|\le
B_{m}\left(v\right)
||R_m||_{\mathbb{T}},\qquad v\in\mathbb{T}.
\]
\end{prop}

Using the relations in Proposition \ref{prop:green_deriv_unitdisk}, 
we can rewrite Proposition \ref{bepr} as follows, where there is no restriction
on the degree of the numerator polynomial in the rational function
(see \cite[Theorem 4]{NK}).

\begin{prop}
\label{thm:genBorweinErdelyi}
Let $R_m\left(v\right)=P\left(v\right)/Q\left(v\right)$
be an arbitrary rational function with no poles on the unit circle,
where $P$ and $Q$ are polynomials.
Denote the poles of $R_m$ by $a_{1},\dots,a_{m} $, where
each pole is repeated as many times as its order.
Then, for $v\in\mathbb{T}$,
\begin{equation*}
|R_m'\left(v\right)|
\le
||R_m||_{\mathbb{T}} 
\cdot\max\left(\sum_{|a_{j}|>1}
\frac{\partial g_{\mathbb{D}_+}\left(v,a_{j}\right)}{\partial
\mathbf{n}_+},
\
\sum_{|a_{j}|<1}
\frac{\partial g_{\mathbb{D}}\left(v,a_{j}\right)}{\partial \mathbf{n}_-}\right).
\end{equation*}
\end{prop}

\subsection{A Gonchar-Grigorjan type estimate}
It is a standard fact that a meromorphic function on a domain with finitely many poles
 can be decomposed into the sum of an analytic function and a rational function
 (which is the sum of the principal parts at the poles). If the rational function
 is required to vanish at $\i$, then this decomposition is unique.

L.D. Grigorjan with A.A. Gonchar
investigated in a series of papers
the supremum norm of the sum of the principal parts of a meromorphic function on the boundary of the given domain
in terms of the supremum norm of the function itself.
In particular, Grigorjan showed in \cite{Grigorjan} that if  $K\subset\mathbf{D}$
is a fixed compact subset of the unit disk $\mathbf{D}$, then
 there exists a constant $C>0$ such that
  all meromorphic functions $f$ on $\mathbf{D}$ having poles only in $K$  have
  principal part $R$ (with $R(\infty)=0$) for which
$\Vert R\Vert\le C \log n \Vert f\Vert$,
where $n$ is the sum of the order of the poles of $f$
 (here $\Vert f\Vert:=\limsup_{|\zeta|\rightarrow 1-} |f(\zeta)|$).

The following recent result (which is \cite[Theorem 1]{multiplygg}) generalizes this to more general domains.
\begin{prop}\label{propgrig}
Suppose that  $D\subset\C$ is a bounded finitely connected domain such that its
boundary $\partial D$ consists of finitely many disjoint $C^{2}$ smooth
Jordan curves. Let $Z\subset D$ be a closed set, and suppose that  $f:D\rightarrow\C$
is a meromorphic function on $D$ such that all of its poles are in $Z$.
Denote the total order of the poles of $f$ by $n$.
If $f_{r}$ is the sum of the principal parts of $f$ (with $f_{r}\left(\infty\right)=0$)
and $f_a$ is its analytic part (so that  $f=f_{r}+f_a$), then
\begin{equation*}
\left\Vert f_{r}\right\Vert _{\partial D},\,\left\Vert f_a\right\Vert _{\partial D}\le C\log n\left\Vert f\right\Vert _{\partial D},
\end{equation*}
where the constant $C=C\left(D,Z\right)>0$ depends only  on $D$ and $Z$.\end{prop}

In this statement
\[\Vert f\Vert_{\partial D}:=\limsup_{\z\in D, \ \zeta\to \partial D} |f(\zeta)|,\]
but we shall apply the proposition in cases when $f$ is actually continuous
on $\partial D$.

\subsection{A Bernstein-Walsh-type approximation theorem}
We shall use the following approximation
theorem.

\begin{prop}\label{propappr} Let
$\tau$ be a Jordan curve and $K$ a compact subset of its interior domain.
Then there are a $C>0$ and $0<q<1$ with
the following property. If $f$ is analytic inside $\tau$
such that $|f(z)|\le M$ for all $z$, then for every
$w_0\in K$ and $m=1,2,\ldots$ there are polynomials $S_m$ of degree
at most $m$ such that $S_m(w_0)=f(w_0)$, $S_m'(w_0)=f'(w_0)$ and
\be \|f-S_m\|_K\le CMq^m.\label{o1}\ee
\end{prop}

\proof. Let $\tau_1$ be a lemniscate, i.e. the level curve of a polynomial,
say $\tau_1=\{z\sep |T_N(z)|=1\}$, such that $\tau_1$ lies inside
$\tau$ and $K$ lies inside $\tau_1$. According to Hilbert's lemniscate theorem
(see e.g. \cite[Theorem 5.5.8]{Ransford})
there is such a $\tau_1$.  Then $K$ is contained in the interior domain of
$\tau_\t=\{z\sep |T_N(z)|=\t\}$ for some $\t<1$. By
Theorem 3 in \cite[Sec. 3.3]{Walsh} (or use \cite[Theorem 6.3.1]{Ransford})
there are polynomials $R_m$ of degree at most $m=1,2,\ldots$ such that
\be  \|f-R_m\|_{\tau_\t}\le C_1Mq^m\label{oop}\ee
with some $C_1$ and $q<1$ (the $q$ depends only on $\t$ and
the degree $N$ of $T_N$). Actually, in that theorem
the right hand side does not show $M$ explicitly, but
the proof, in particular the error formula (12) in \cite[Section 3.3]{Walsh}
(or the error formula (6.9) in \cite[Section 6.3]{Ransford}),  gives
(\ref{oop}).

Now (\ref{oop})  pertains to hold
also on the interior domain to $\tau_\t$, so if $\d$  is the
distance in between $K$ and $\tau_\theta$ and $w_0\in K$, then
for all $|\xi-w_0|=\d$ we have $|f(\xi)-R_m(\xi)|\le C_1Mq^m$.
Hence, by Cauchy's integral formula for the derivative we have
\[|f'(w_0)-R_m'(w_0)|\le \frac{C_1Mq^m}{\d}.\]
Therefore, the polynomial
\[S_m(z)=R_m(z)+(f(w_0)-R_m(w_0))+(f'(w_0)-R_m'(w_0))(z-w_0)\]
satisfies the requirements with $C=C_1(2+{\rm diam}(K)/\d)$ in (\ref{o1}).\endproof

\subsection{Bounds and smoothness for Green's functions}\label{sectboundgreen}

In this section we collect some simple facts on Green's functions and their normal
derivatives.

Let $K\subset\C$ be a compact set with connected complement and $Z\subset \C\setm K$
a closed set. Suppose that $\s$ is a Jordan curve that separates $K$ and $Z$, say
$K$ lies in the interior of $\s$ while $Z$ lies in its exterior.
Assume also that there is a family $\{\g_\tau\}\subset K$ of Jordan arcs
such that ${\rm diam}(\g_\tau)\ge d>0$ with some $d>0$, where
${\rm diam}(\g_\tau)$ denotes the diameter of $\g_\tau$.

First we prove
\begin{prop}\label{lemgreencomp} There are $c_0,C_0>0$ such that for
all $\tau$, $z\in \s$ and all $a\in Z$ we have
\be c_0\le g_{\C\setm \g_\tau}(z,a)\le C_0.\label{c1}\ee
\end{prop}
\proof. We have the formula (\cite[p. 107]{Ransford})
\[g_{\C\setm \g_\tau}(z,\i)=\log\frac{1}{\c{(\g_\tau)}}+\int\log|z-t|d\mu_{\g_\tau}(t),\]
where $\mu_{\g_\tau}$ is the equilibrium measure of $\g_\tau$ and where $\c(\g_\tau)$ denotes
the logarithmic capacity of $\g_\tau$. Since (see \cite[Theorem 5.3.2]{Ransford})
\[\c(\g_\tau)\ge \frac{{\rm diam}(\g_\tau)}{4}\ge \frac{d}{8},\]
and for $z\in \s$, $t\in \g_\tau$ we have $|z-t|\le {\rm diam}(\s)$, we obtain
\[g_{\C\setm \g_\tau}(z,\i)\le \log\frac{1}{d/8}+\log {\rm diam}(\s)=:C_1.\]

Let $\O$ be the exterior of $\s$ (including $\i$). By Harnack's inequality (\cite[Corollary 1.3.3]{Ransford})
for any closed set
$Z\subset \O$ there is a constant $C_Z$ such that for all positive harmonic functions
$u$ on $\O$ we have
\[\frac{1}{C_Z}u(\i)\le u(a)\le C_Zu(\i),\qquad a\in Z.\]
Apply this to the harmonic function $g_{\C\setm \g_\tau}(z,a)=
g_{\C\setm \g_\tau}(a,z)$ (recall that Green's functions are symmetric in their
arguments), $z\in \s$, $a\in Z$, to conclude for $z\in \s$
\[g_{\C\setm \g_\tau}(z,a)=g_{\C\setm \g_\tau}(a,z)\le C_Z g_{\C\setm \g_\tau}(\i,z)
=C_Z g_{\C\setm \g_\tau}(z,\i)\le C_ZC_1.\]

To prove a lower bound note that
\[g_{\C\setm \g_\tau}(z,\i)\ge
 g_{\C\setm K}(z,\i)\ge c_1,\qquad z\in \s,\]
 because $\g_\tau\subset K$ and  $g_{\C\setm K}(z,\i)$ is a positive
 harmonic function outside $K$. From here we get
\[g_{\C\setm \g_\tau}(z,a)\ge
\frac{c_1}{C_Z},\qquad z\in \s,\ a\in Z,\]
exactly as before by appealing to the symmetry of
the Green's function and to Harnack's inequality.\endproof

\begin{cor} \label{corgreencomp} With the $c_0,C_0$ from the preceding
lemma for
all $\tau$, $a\in Z$ and for all $z$ lying inside
$\s$ we have
\be \frac{c_0}{C_0}g_{\C\setm \g_\tau}(z,\i)\le g_{\C\setm \g_\tau}(z,a)\le
\frac{C_0}{c_0}g_{\C\setm \g_\tau}(z,\i).\label{c2}\ee
\end{cor}
\proof. For $z\in \s$ the inequality (\ref{c2}) was shown
in the preceding proof. Since both
$g_{\C\setm \g_\tau}(z,\i)$ and $g_{\C\setm \g_\tau}(z,a)$
are harmonic in the domain that lies in between $\g_\tau$ and $\s$
and both vanish on $\g_\tau$, the statement follows
from the maximum principle.\endproof

Next, let $\G$ be a $C^2$ Jordan curve and $G_\pm$ the interior and exterior domains to $\G$ (see Section \ref{sectresults}).
Assume, as before, that $Z\subset\C\setm \G$ is a closed set.

\begin{prop}\label{propgreencont} There are constants $C_1, c_1>0$ such that
\be c_1\le \frac{\partial g_{G_-}(z_0,a)}{\partial {\bf n}_-}\le C_1,\qquad a\in Z\cap G_-
\label{b41}\ee
and
\be c_1\le \frac{\partial g_{G_+}(z_0,a)}{\partial {\bf n}_+}\le C_1,\qquad a\in Z\cap G_+.
\label{b410}\ee
These bounds hold uniformly in $z_0\in \G$. Furthermore, the Green's functions $g_{G_\pm}(z,a)$,
$a\in Z$,  are uniformly H\"older 1 equicontinuous close to the boundary $\G$.
\end{prop}

\proof. It is enough to prove (\ref{b41}). Let $b_0\in G_-$ be a fixed point and let $\f$ be a conformal map from
the unit disk ${\mathbb{D}}$ onto $G_-$ such that $\f(0)=b_0$. By the Kellogg-Warschawski theorem
(see \cite[Theorem 3.6]{Pommerenke})
$\f'$ has a continuous extension to the closed unit disk which does not vanish there.
It is clear that $g_{G_-}(z,b_0)=
-\log|\f^{-1}(z)|$, and consider some local branch of $-\log \f^{-1}(z)$
for $z$ lying close to $z_0$. By the
Cauchy-Riemann equations
\[\frac{\partial g_{G_-}(z_0,b_0)}{\partial {\bf n}_-}=\left|{\left(-\log\f^{-1}(z)\right)'\restr{z=z_0}}\right|\]
(note that the directional derivative of $g_{G_-}$ in the direction perpendicular
to ${\bf n}_-$ has 0 limit at $z_0\in \partial G_-$),
so we get the formula
\be \frac{\partial g_{G_-}(z_0,b_0)}{\partial {\bf n}_-}=\frac{1}{|\f'(\f^{-1}(z_0)|},\label{37}\ee
which shows that this normal derivative is finite, continuous in $z_0\in \G$ and positive.

Let now $\s$ be a Jordan curve that separates $(Z\cap G_-)\cup\{b_0\}$ from $\G$. Map $G_-$ conformally onto $\C\setm [-1,1]$
by a conformal map $\F$ so that $\F(b_0)=\i$. Then $g_{G_-}(z,a)=g_{\C\setm [-1,1]}(\F(z),\F(a))$,
and $\F(\s)$ is a Jordan curve that separates $\F((Z\cap G_-)\cup\{b_0\})$
from $[-1,1]$. Now apply Proposition \ref{lemgreencomp} to $\C\setm [-1,1]$ and to $\F(\s)$ to conclude that
all the Green's functions $g_{\C\setm [-1,1]}(w,\F(a))$,  $a\in Z\cup \{b_0\}$,
are comparable
on $\F(\s)$ in the sense that all of them lie in between two positive constants $c_2<C_2$
there. In view of what we have just said, this means that the Green's functions
$g_{G_-}(z,a)$, $a\in Z\cup \{b_0\}$, are comparable
on $\s$ in the sense that all of them lie in between the same  $c_2<C_2$
there. But then, as in Corollary \ref{corgreencomp}, they are also comparable in the domain that lies in
between $\G$ and $\s$, and hence
\[\frac{c_2}{C_2}\frac{\partial g_{G_-}(z_0,b_0)}{\partial {\bf n}_-}\le
\frac{\partial g_{G_-}(z_0,a)}{\partial {\bf n}_-}
\le \frac{C_2}{c_2}\frac{\partial g_{G_-}(z_0,b_0)}{\partial {\bf n}_-},\qquad a\in Z,\]
which proves (\ref{b41}) in view of (\ref{37}).

The uniform H\"older continuity is also easy to deduce
from (\ref{37}) if we compose $\f$ by fractional linear mappings
of the unit disk onto itself (to move the pole $\f(0)$ to other points).\endproof

\sect{The Bernstein-type inequality on analytic curves}\label{sectanalcurve}
\label{sectproofmaincurve}

In this section we assume that $\G$ is analytic, and prove
(\ref{maincurve}) using Propositions \ref{thm:genBorweinErdelyi}, \ref{propgrig} and \ref{propappr}.

Fix $z_0\in \G$ and consider the conformal maps $\F_1$ and $\F_2$ from
Section \ref{sectconf}. Recall that the inner map $\F_1$ has an extension
to a disk ${D_1}=\{z\sep |z|< r_1\}$ and
the external map $\F_2$ has an extension
to the exterior ${D_2}=\{z\sep |z|> r_2\}$ of a disk
with some $r_2<1<r_1$. For simpler notation, in what follows
we shall assume that $\F_1$ resp. $\F_2$ actually have extensions to
a neighborhood of the closures $\ov{D_1}$ resp. $\ov{D_2}$ (which
can be achieved by decreasing $r_1$ and increasing $r_2$ if necessary).

In what follows we set ${\mathbb{T}}(r)=\{z\sep |z|=r\}$ for the circle
of radius $r$ about the origin. As before, ${\mathbb{T}}={\mathbb{T}}(1)$ denotes the
unit circle.

The constants $C,c$ below depend only  on $\Gamma$ and they are not
the same at each occurrence.

\smallskip

We decompose $R_n$ as,
\[
R_n=f_{1}+f_{2}
\]
where $f_{1}$ is a rational function with poles in $Z\cap G_-$,
$f_{1}\left(\infty\right)=0$
and $f_{2}$ is a rational function with poles in $Z\cap G_+$.
This decomposition is unique.
If we put $N_1:=\deg\left(f_1\right)$,
$N_2:=\deg\left(f_2\right)$, then $N_1+N_2=n$.
Denote the poles of $f_1$ by
$\alpha_{j}$, $j=1,\ldots,N_1$,
and
the poles of $f_2$ by
$\beta_{j}$, $j=1,\ldots,N_2$
 (with counting the orders of the poles).

We use Proposition \ref{propgrig}
on $G_-$
to conclude
\begin{equation}
\left\Vert f_{1}\right\Vert_{\Gamma},
\left\Vert f_{2}\right\Vert_{\Gamma}
\le
C 
\log n
\left\Vert R_n\right\Vert_{\Gamma}.
\label{eq:f_onetwo_norm_est}
\end{equation}
By the maximum modulus principle then it follows that
\begin{equation}
\left\Vert f_1\right\Vert_{\outercurve}
\leq C 
\log n
\left\Vert R_n\right\Vert_\Gamma
\label{est:f_one_on_outer_curve}
\end{equation}
and
\begin{equation}
\left\Vert f_2\right\Vert_{\innercurve}
\leq C  
\log n
\left\Vert R_n\right\Vert_\Gamma.
\label{est:f_two_on_inner_curve}
\end{equation}

Set $F_1:=f_1(\Phi_1)$ and $F_2:=f_2(\Phi_2)$.
These are meromorphic functions in $\innerdisk$ and $\outerdisk$ resp.
with poles at $\alpha_j':=\Phi_1^{-1}(\alpha_j)$, $j=1,\ldots,N_1$ and
at $\beta_k':=\Phi_2^{-1}(\beta_k)$, $k=1,\ldots,N_2$.

Let $F_1=F_{1,r}+F_{1,a}$ be the decomposition of $F_1$ with
respect to the unit disk into rational and analytic parts with
$F_{1,r}(\infty)=0$, and in a similar
fashion, let $F_2=F_{2,r}+F_{2,a}$ be the decomposition of $F_2$ with
respect to the exterior of the unit disk into rational and analytic parts with
$F_{2,r}(0)=0$.
(Here $r$ refers to the rational part, $a$ refers to the analytic part.)
Hence, we have by Proposition \ref{propgrig}
\[\left\Vert F_{j,r}\right\Vert_{\unitcircle}, \left\Vert
F_{j,a}\right\Vert_{\unitcircle} \leq C 
\log n
\left\Vert F_{j}\right\Vert_{\unitcircle}, \qquad j=1,2.\]
Thus, $F_{1,r}$ is a rational function with poles at $\alpha_j'\in {\mathbb{D}}$,
so by the maximum modulus theorem and \eqref{eq:f_onetwo_norm_est} 
we have
\begin{equation}
\left\Vert F_{1,r} \right\Vert_{\outercirc} \leq \left\Vert F_{1,r}\right\Vert_{\unitcircle}\leq
C 
\log n  \left\Vert F_1\right\Vert_{\unitcircle}
\leq
C 
\log^2 n  \left\Vert R_n\right\Vert_{\G},
\label{est:foner_norm}
\end{equation}
where we used that $\|F_1\|_{\mathbb{T}}=\|f_1\|_\G$.
But \eqref{est:f_one_on_outer_curve}
is the same as
\begin{equation*}
\left\Vert F_1 \right\Vert_{\outercirc}
\leq
C 
\log n  \Vert R_n\Vert_\Gamma,
\end{equation*}
so we can conclude also
\begin{equation}
\left\Vert F_{1,a} \right\Vert_{\outercirc}
\leq
C 
\log^2 n  \left\Vert R_n\right\Vert_{\Gamma}.
\label{est:fonea_norm}
\end{equation}

Thus, $F_{1,a}$ is an analytic function in $\innerdisk$ with the bound
in \eqref{est:fonea_norm}. Apply now Proposition \ref{propappr} to this function
and to the unit circle as $K$ (and with a somewhat larger concentric circle as
$\tau$) with degree $m=[\sqrt n]$. According to that proposition there are $C,c>0$
and polynomials $S_1=S_{1,\sqrt n}$ of degree at most $\sqrt n$ such that
\begin{equation*}
\left\Vert F_{1,a}-S_1 \right\Vert_{\unitcircle}
\leq C e^{-c \sqrt{n}}
\Vert R_n \Vert_\Gamma,
\quad
S_1(1)=F_{1,a}(1),
\
S_1'(1)=F_{1,a}'(1).
\end{equation*}
Therefore, $\tilde {\cal R}_1:=F_{1,r}+S_1$ is a rational function with poles at
$\alpha_j'$, $j=1,\ldots,N_1$
and with a pole at $\infty$ with order at most $\sqrt{n}$ which satisfies
\begin{equation}
\left\Vert F_1 -\tilde {\cal R}_1 \right\Vert_{\unitcircle}
\leq C e^{-c\sqrt{n}} \Vert R_n \Vert_\Gamma, \quad
\tilde {\cal R}_1(1)=F_1(1),
\
\tilde {\cal R}_1'(1) = F_1'(1)
\label{rel:fone_rone}
\end{equation}

In a similar vein, if we consider $F_2(1/v)$ and use
\eqref{est:f_two_on_inner_curve},
then we get a polynomial $S_2$
of degree at most $\sqrt{n}$ such that
\begin{equation*}
\left\Vert F_{2,a}(1/v)-S_2(v) \right\Vert_{\unitcircle}
\leq
C e^{-c\sqrt{n}} \Vert R_n \Vert_\Gamma, \quad
S_2(1)=F_{2,a}(1),
\
S_2'(1)=- F_{2,a}'(1)
\end{equation*}
But then $\tilde {\cal R}_2(v):=F_{2,r}(v)+S_2(1/v)$ is a rational function
with poles at $\beta_k'$, $k=1,\ldots, N_2$ and with a pole at $0$ of
order at most $\sqrt{n}$
that satisfies
\begin{equation}
\left\Vert F_2 - \tilde {\cal R}_2 \right\Vert_{\unitcircle}\leq
C e^{-c\sqrt{n}} \Vert R_n\Vert_\Gamma, \quad
\tilde {\cal R}_2(1)=F_2(1),\
\tilde {\cal R}_2'(1)=F_2'(1).
\label{rel:ftwo_rtwo}
\end{equation}

What we have obtained is that the rational function $\tilde {\cal R}:=\tilde {\cal R}_1+\tilde {\cal R}_2$
is of distance $\leq C e^{-c\sqrt{n}}\Vert R_n \Vert_\Gamma$
from $F_1+F_2$ on the unit circle
and it satisfies
\begin{equation}
\tilde {\cal R}(1)=\left(F_1+F_2\right)(1)
= f_1(z_0)+f_2(z_0)=R_n (z_0)\label{480}
\end{equation}
and using \eqref{eq:phione_phitwo},
\begin{equation}
\tilde {\cal R}'(1)=
\left(F_1'+F_2'\right)(1)
= f_1'(z_0)\Phi_1'(1)+f_2'(z_0)\Phi_2'(1)
=R_n'(z_0)\Phi_1'(1).
\label{eq:roneprime}
\end{equation}

Consider now $F_1+F_2$ on the unit circle, i.e.
\begin{equation*}
F_1(e^{it})+F_2(e^{it})=
f_1(\Phi_2(e^{it}))+
f_2(\Phi_2(e^{it}))+
f_1(\Phi_1(e^{it}))-f_1(\Phi_2(e^{it})).
\end{equation*}
The sum of the first two terms on the right is $R_n(\Phi_2(e^{it}))$,
and this is at most $\Vert R_n\Vert_\Gamma$ in absolute value.
Next, we estimate the difference of the last two terms.

The function $\Phi_1(v)-\Phi_2(v)$
is analytic in the ring  $r_2<|v|<r_1$ and it is bounded there
with a bound depending only on $\Gamma,r_1,r_2$,
furthermore it has a double zero at $v=1$ (because
of (\ref{eq:phione_phitwo})).  These imply
\[|\Phi_1(e^{it})-\Phi_2(e^{it})|\le C|e^{it}-1|^2\le Ct^2, \qquad
t\in [-\pi,\pi],\]
with some constant $C$. By Proposition \ref{prop:rough_berns}
we have with (\ref{eq:f_onetwo_norm_est}) also
the bound
\[\|f_1'\|_{\Gamma}\le Cn\log n \|R_n\|_{\Gamma},\]
and these last two facts give us (just integrate $f_1'$ along the shorter
arc of $\G$ in between $\Phi_1(e^{it})$ and $\Phi_2(e^{it})$ and use that the
length of this
arc is at most $C|\Phi_1(e^{it})-\Phi_2(e^{it})|$)
\[|f_1(\Phi_1(e^{it}))-f_1(\Phi_2(e^{it}))|\le Ct^2n\log n \|R_n\|_{\Gamma}.\]

By \cite[Theorem 4.1]{TotikTrans} there are polynomials $Q$ of degree
at most $[n^{4/5}]$ such that
 $Q(1)=1$, $\|Q\|_{\mathbb{T}}\le 1$,  and with some constants $c_0,C_0>0$
\[
|Q(v)|\le C_0\exp(-c_0 n^{4/5}|v-1|^{3/2}),\qquad |v|=1.
\]
With this $Q$ consider the rational function
${\cal R}(v)=\tilde {\cal R}(v) Q(v)$. On the unit circle this is closer
than  $Ce^{-c\sqrt n}\|R_n\|_{\Gamma}$ to $(F_1+F_2)Q$, and in view of
what we have just proven,
we have at $v=e^{it}$
\[
|(F_1(v)+F_2(v))Q(v)|\le \|R_n\|_{\Gamma}+ Ct^2n\log n
 C_0\exp\left(-c_0 n^{4/5}|t/2|^{3/2}\right)\|R_n\|_{\Gamma}.
\]
On the right
\begin{multline*}
t^2n\log n \exp\left(-c_0 n^{4/5}|t/2|^{3/2}\right)
\\
= 4\left(n^{4/5}|t/2|^{3/2}\right)^{4/3}\exp\left(-c_0 n^{4/5}|t/2|^{3/2}\right)
\frac{\log n}{n^{1/15}}
\le  C\frac{\log n}{n^{1/15}}
\end{multline*}
because $|x|^{4/3}\exp(-c_0 |x|)$ is bounded on the real line.

All in all, we obtain
\be
\|{\cal R}\|_{\unitcircle}\le (1+o(1))\|R_n\|_{\Gamma},\label{iii}
\ee
and
\begin{multline*}
|{\cal R}'(1)|=|\tilde {\cal R}'(1)Q(1)+\tilde {\cal R}(1)Q'(1)|
=|\tilde {\cal R}'(1)|+O\left(|\tilde {\cal R}(1)||Q'(1)|\right)
\\
=|R_n'(z_0)|+O(n^{4/5})\|R_n\|_\G,
\end{multline*}
where we used $Q(1)=1$, (\ref{480})--(\ref{eq:roneprime}), $|\Phi_1'(1)|=1$
and the classical Bernstein inequality
for $Q'(1)$, which gives the bound $n^{4/5}$ for the derivative of $Q$.

The poles of ${\cal R}$ are at $\alpha_j'$, $1\le j\le N_1$, and at
$\beta_k'$, $1\le k\le N_2$,
as well as a $\le n^{1/2}$ order pole at 0 (coming from
the construction of $S_{2,n}$) and a $\le n^{1/2}+n^{4/5}$
order pole at $\infty$ (coming from the construction of $S_{1,n}$
and the use of $Q$).

Now we apply the Borwein-Erd\'elyi inequality (Proposition
\ref{thm:genBorweinErdelyi})
to $|{\cal R}'(1)|$ to obtain
\bean |R_n'(z_0)|&\le& |{\cal R}'(1)|+O(n^{4/5})\|R_n\|_\G\\
&\le &\|{\cal R}\|_{\mathbb{T}}\max\left(\sum_{k}\frac{\partial g_{{\mathbb{D}}_+}(1,\b_k')}{\partial {\bf n}_+}
+(n^{1/2}+n^{4/5})\frac{\partial g_{{\mathbb{D}}_+}(1,\i)}{\partial {\bf n}_+},\right.\\
&&\left.\sum_{j}\frac{\partial g_{{\mathbb{D}}}(1,\a_j')}{\partial {\bf n}_-}
+n^{1/2}\frac{\partial g_{{\mathbb{D}}}(1,0)}{\partial {\bf n}_-}\right)
+O(n^{4/5})\|R_n\|_\G.\eean
If we use here how the normal derivatives
transform under the mappings $\Phi_1$ and $\Phi_2$ as in
Proposition  \ref{prop:green_deriv_unitdisk},
then we get from (\ref{iii})
\bean |R_n'(z_0)|&\le& (1+o(1))\|R_n\|_{\G}\times\\[10pt]
&&\max\left(\sum_{a\in Z\cap G_+}\frac{\partial g_{G_+}(z_0,a)}{\partial {\bf n}_+}
+(n^{1/2}+n^{4/5})\frac{\partial g_{G_+}(z_0,\F_2(\i))}{\partial {\bf n}_+},\right.\\[10pt]
&&\left.\sum_{a\in Z\cap G_-}\frac{\partial g_{G_-}(z_0,a)}{\partial {\bf n}_-}
+n^{1/2}\frac{\partial g_{G_-}(z_0,\F_1(0))}{\partial {\bf n}_-}\right)
+O(n^{4/5})\|R_n\|_\G.\eean
Since, by (\ref{b41})--(\ref{b410}), the normal derivatives on the right lie
in between two positive constants that depend only on $\G$ and $Z$,
(\ref{maincurve}) follows (note that one of the sums
$\sum_{a\in Z\cap G_+}$ or $\sum_{a\in Z\cap G_-}$ contains at least $n/2$ terms).\endproof

\sect{The Bernstein-type inequality on analytic arcs}\label{sectanalarcs}
In this section we prove Theorem \ref{tharc} in the case when
the arc $\G$ is analytic. We shall reduce this case to Theorem \ref{thcurve}
for analytic Jordan curves that has been proven in the preceding section.
We shall use the Joukowskii map
to transform the arc setting  to the curve setting.

For clearer notation let us write for the arc in Theorem \ref{tharc} $\G_0$.
We may assume that the endpoints of $\G_0$ are $\pm 1$. Consider the
pre-image $\G$ of $\G_0$ under the Joukowskii map
$z=F(u)=(u+1/u)/2$. Then $\G$ is a Jordan curve, and if
$G_\pm$ denote the inner and outer domains to $\G$, then
$F$ is a conformal map from both $G_-$ and from $G_+$ onto
$\C\setm \G_0$. Furthermore, the analyticity of $\G_0$ implies
that $\G$ is an analytic Jordan curve, 
this follows from standard steps, see e.g. \cite{NK}, p. 875.

Denote the inverse of $z=F(u)$ restricted to $G_-$ by
$F_1^{-1}(z)=u$ and
that restricted to $G_+$
by $F_2^{-1}(z)=u$. So $F_j(z)=z\pm \sqrt{z^2-1}$
with an appropriate branch of $\sqrt{z^2-1}$ on the plane
cut along $\G_0$.

\begin{figure}
\begin{center}
\includegraphics[keepaspectratio,width=\textwidth]{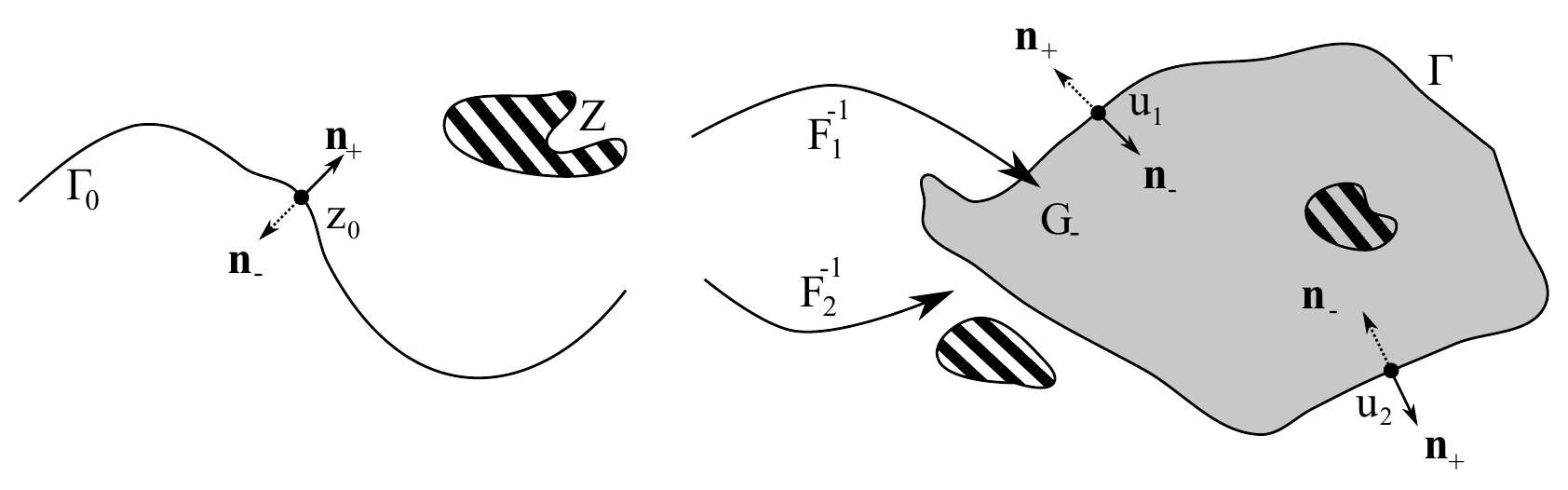}
\end{center}
\caption{The open-up}
\label{fig:openup}
\end{figure}

We need the mapping properties of $F$
regarding normal vectors,
for full details, we refer to \cite{NK} p. 879.
Briefly, for any $z_0\in \G_0$ that is not one of the endpoints
of $\G_0$ there are exactly two $u_1,u_2\in \Gamma$, $u_1\ne u_2$
such that $F(u_1)=F(u_2)=z_0$.
Denote the normal vectors to $\Gamma$
pointing outward by $\mathbf{n}_+$ and the normal vectors pointing
inward by $\mathbf{n}_-$ (it is usually unambiguous from the context
at which point $u\in\Gamma$ we are referring to).
By reindexing  $u_{1}$ and $u_{2}$ (and possibly reversing the
parametrization of $\Gamma_0$),
we may assume that the  (direction of the) normal
vector $\mathbf{n}_+(u_1)$ is mapped by $F$ to the  (direction of the) normal
vector $\mathbf{n}_+(z_0)$.
This then implies that  (the directions of) $\mathbf{n}_+(u_1)$,
$\mathbf{n}_-(u_1)$  and $\mathbf{n}_+(u_2)$, $\mathbf{n}_-(u_2)$
 are mapped
by $F$ to  (the directions of) $\mathbf{n}_+$,
$\mathbf{n}_-$,
$\mathbf{n}_-$,
$\mathbf{n}_+$ at $z_0$,
respectively.
These mappings are depicted on Figure \ref{fig:openup}.

The corresponding normal derivatives of the Green's
functions are related as follows.
\begin{prop}
\label{prp:greenopenup}
We have for $a\in \C\setm \G$
\begin{multline*}
\frac{\partial g_{\C\setminus \Gamma_0}\left(z_0,a\right)}{\partial \mathbf{n}_-}
=
\frac{\partial g_{G_{-}}\left(u_{1},F_{1}^{-1}\left(a\right)\right)}{\partial \mathbf{n}_-}
/\left|F'\left(u_{1}\right)\right|
\\ =
\frac{\partial g_{G_{+}}\left(u_{2},F_{2}^{-1}\left(a\right)\right)}{\partial \mathbf{n}_+}
/\left|F'\left(u_{2}\right)\right|
\end{multline*}
and, similarly for the other side,
\begin{multline*}
\frac{\partial g_{\C\setminus \Gamma_0}\left(z_0,a\right)}{\partial \mathbf{n}_+}
=
\frac{\partial g_{G_{-}}\left(u_{2},F_{1}^{-1}\left(a\right)\right)}{\partial \mathbf{n}_-}
/\left|F'\left(u_{2}\right)\right|
\\ =
\frac{\partial g_{G_{+}}\left(u_{1},F_{2}^{-1}\left(a\right)\right)}{\partial \mathbf{n}_+}
/\left|F'\left(u_{1}\right)\right|.
\end{multline*}
\end{prop}
This proposition follows immediately from
\cite[Proposition 6]{NK} and is a two-to-one mapping analogue of Proposition
\ref{prop:green_deriv_unitdisk}.
\medskip

After these preliminaries let us turn to the proof of (\ref{mainarc})
at a point $z_0\in \G_0$. Consider $f_1(u):=R_n(F(u))$ on the analytic Jordan curve $\Gamma$
at $u_1$ (where $F(u_1)=z_0$). This is a rational function
with poles at $F_1^{-1}(a)\in G_-$ and at
$F_2^{-1}(a)\in G_+$, where $a$ runs through the poles of $R_n$.
According to (\ref{maincurve}) (that has been
verified in Section \ref{sectanalcurve} for analytic curves)
we have
\begin{multline*}
\left|f_1' \left(u_1\right)\right|
\le
(1+o(1))
\left\|f_1\right\|_{\Gamma}
\\ \cdot
\max\left(
\sum_a
\frac{\partial g_{G_-}\left(u_1, F_1^{-1}(a)\right)}{\partial \mathbf{n}_-},
\sum_a
\frac{\partial g_{G_+}\left(u_1, F_2^{-1}(a)\right)}{\partial \mathbf{n}_+}
\right),
\end{multline*}
where $a$ runs through the poles of $R_n$ (counting multiplicities).
If we use here that
$\left\|f_1\right\|_{\Gamma}=\|R_n\|_{\G_0}$ and
$f_1'\left(u_1\right)=R_n'\left(z_0\right) F'(u_1)$,
we get from  Proposition \ref{prp:greenopenup}
\begin{multline*}
\left|R_n' \left(z_0\right)\right|
\le
(1+o(1))
\left\|R_n\right\|_{\Gamma_0}
\\ \cdot
\max\left(
\sum_a
\frac{\partial g_{\C\setm \G_0}\left(z_0, a\right)}{\partial \mathbf{n}_-},
\sum_a
\frac{\partial g_{\C\setm \G_0}\left(z_0, a\right)}{\partial \mathbf{n}_+}
\right),
\end{multline*}
which is (\ref{mainarc}) when $\G$ is replaced by $\G_0$.\endproof

\sect{Proof of Theorem \ref{tharc}}\label{sectc2arcs}
In this section we verify (\ref{mainarc}) for $C^2$ arcs. Recall that
in Section \ref{sectanalarcs}  (\ref{mainarc}) has already been proven for
analytic arcs and we shall reduce the $C^2$ case to that by approximation
similar to what was used in \cite{TotikMarkov}.

In the proof we shall frequently identify a Jordan arc
with its parametric representation.

By assumption,
$\G$ has a twice differentiable
parametrization $\g(t)$, $t\in[-1,1]$, such that
$\g'(t)\not=0$ and $\g''$ is continuous.
We may assume
that $z_0=0$ and that the real line is tangent to $\G$ at $0$, and
also
that $\g(0)=0$, $\g'(0)>0$. There is an $M_1$ such that
for all $t\in [-1,1]$
\be \frac1{M_1}\le |\g'(t)|\le M_1,\qquad |\g''(t)|\le M_1.\label{0}\ee

Let $\g_0:=\g$, and for some $0<\tau_0<1$ and for all $0<\tau\le \tau_0$
choose a polynomial $g_\tau$ such that
\be |\g''-g_\tau|\le \tau,\label{alt0}\ee
and set
\be \g_\tau(t)=\int_0^t\left(\int_0^u g_\tau(v)dv+\g_0'(0)\right) du.
\label{gdef}\ee
It is clear that
\be |\g_\tau(t)-\g_0(t)|\le \tau |t|^{2},\quad |\g_\tau'(t)-\g_0'(t)|\le \tau |t|,
\quad |\g_\tau''(t)-\g_0''(t)|\le \tau.\label{alt1}\ee

It was proved in \cite[Section 2]{TotikMarkov} that for small $\tau$,
say for all $\tau\le \tau_0$ (which can be achieved by decreasing $\tau_0$
if necessary), these $\g_\tau$ are analytic Jordan arcs, and
\be  g_{\C\setm \g_0}(z,\i)\le M_2\sqrt\tau |z|^2,
\quad z\in \g_\tau,
\label{v2**}\ee
with some constant $M_2$ that is independent of $\tau$ and $z$.
We need similar estimates for all $g_{\C\setm \g_\tau}(z,a)$, $a\in Z$.
To get them consider the closure of the set
$\cup_{0\le \tau\le \tau_0}\g_\tau$ and its polynomial convex hull
\[K={\rm Pc}\left(\overline{\bigcup_{0\le \tau\le \tau_0}\g_\tau}\right),\]
which is the union of that closure with all the bounded components
of its complement. Now this is a situation when
the results from Section \ref{sectboundgreen} can be applied.
From Corollary \ref{corgreencomp}
and from (\ref{v2**}) we can conclude for
all $a\in Z$
\be  g_{\C\setm \g_0}(z,a)\le M_3\sqrt\tau |z|^2,
\quad z\in \g_\tau
\label{v2**+}\ee
with some constant $M_3$.

Let ${\bf n}_\pm$ denote the two normals to $\g_\tau$
 at the origin. Note that ${\bf n}_\pm$ are common to all the arcs
 $\g_\tau$, $0\le \tau\le \tau_0$.

\begin{lem}\label{poscor} For small $\tau_0$ the normal derivatives
\[\frac{\partial g_{\C\setm \g_\tau}(0,a)}{\partial {\bf n}_\pm},\qquad 0\le \tau\le \tau_0,\ a\in Z\cup\{\i\},\]
are uniformly bounded from below and above by a positive number.
\end{lem}
\proof. It was proven in \cite[Appendix 1]{TotikMarkov}
that
\be \frac{\partial g_{\C\setm \g_\tau}(0,\i)}{\partial {\bf n}_\pm}\to
\frac{\partial g_{\C\setm \g_0}(0,\i)}{\partial {\bf n}_\pm}\label{111}\ee
as $\tau\to 0$, and the value on the right is positive and finite. Now
just invoke Corollary \ref{corgreencomp} (note that (\ref{c2}) implies similar
inequalities for the normal derivatives).\endproof

Next we mention that  (\ref{alt1}) implies
the following:  no matter how
$\eta>0$ is given, there is a $\tau_\eta<\tau_0$ such that for
$\tau<\tau_\eta$ we have
\be \frac{\partial g_{\C\setm \g_\tau}(0,\i)}{\partial {\bf n}_\pm}<
(1+\eta)\frac{\partial g_{\C\setm \g_0}(0,\i)}{\partial {\bf n}_\pm}.
\label{normal}\ee
 In fact, (\ref{111}) was proven in \cite[Appendix 1, (6.1)]{TotikMarkov}
under the assumption
  (\ref{alt1}), and
since the  normal derivatives on the right are not zero, (\ref{normal})
follows.

We shall also need this inequality when $\i$ is replaced by
an arbitrary pole $a\in Z$.
Let $a\in Z$ be arbitrary, and consider
the mapping $\f_a(z)=1/(z-a)$. Under this mapping
$\g_\tau$ is mapped into $\f_a(\g_\tau)$ with parametrization
$\f_a(\g_\tau(t))$, $t\in [-1,1]$, and it is clear that
(\ref{alt1}) implies its analogue for the image curves:
\bean |\f_a(\g_\tau)(t)-\f_a(\g_0)(t)|\le C\tau |t|^{2},&\quad& |(\f_a(\g_\tau))'(t)-
(\f_a(\g_0))'(t)|\le C\tau |t|,\\
&\quad& |(\f_a(\g_\tau))''(t)-(\f_a(\g_0))''(t)|\le C\tau,\eean
with some constant $C$ that is independent of $\tau$ and $a\in Z$.
Furthermore,
\[g_{\C\setm \g_\tau}(z,a)=g_{\C\setm \f_a(\g_\tau)}(\f_a(z),\i),\]
\[\frac{\partial g_{\C\setm \g_\tau}(0,a)}{\partial {\bf n}_\pm}
=\frac{\partial g_{\C\setm \f_a(\g_\tau)}(\f_a(0),\i)}{\partial {\bf n}(\f_a(0))_\pm}|\f_a'(0)|.\]
Now if we use these in the proof of
\cite[Appendix 1]{TotikMarkov} and use also Lemma \ref{poscor}, then
we obtain that for every
$\eta>0$ there  is a $\tau_\eta<\tau_0$ such that for
$\tau<\tau_\eta$ and $a\in Z$ we have
\be \frac{\partial g_{\C\setm \g_\tau}(0,a)}{\partial {\bf n}_\pm}<
(1+\eta)\frac{\partial g_{\C\setm \g_0}(0,a)}{\partial {\bf n}_\pm}.
\label{normal11}\ee
An inspection of the proof reveals that $\tau_\eta$ can be made independent
of $a\in Z$, so (\ref{normal11}) is uniform in $a\in Z$.

After these preparations let $R_n$ be a rational function with
poles in $Z$ such that total order of its poles
(including possibly the pole at $\i$) is $n$.
We use
\be |R_n(z)|\le \exp\left(\sum_a g_{\C\setm \G}(z,a)\right)\|R_n\|_\G,
\label{aa}\ee
where the summation is for all poles of $R_n$. This is the analogue of (\ref{BW}),
and its proof is the same that we gave for (\ref{BW}).
 Hence, in view of
 (\ref{v2**+}), we
have   for $z\in \g_\tau$ (recall that $\g_0=\G$)
\be |R_n(z)|\le \|R_n\|_{\G} \exp\Bigl(nM_3\sqrt\tau |z|^2\Bigr).
\label{ee}\ee

The polynomial convex hull $K$ introduced above has the property
 that there is a disk (say in the upper half plane)
in the complement of $K$ which contains the point 0 on
its boundary. Indeed, this easily follows
from the  construction of the curves $\g_\tau$.
Now we use \cite[Theorem 4.1]{TotikTrans}, according to which
  there are constants
 $c_1,C_1$ and for each $m$ polynomials
$Q_m$ of degree at most $m$ such that
\begin{description}
\item[(i)] $Q_m(0)=1$,
\item[(ii)] $|Q_m(z)|\le 1,\quad z\in K$,
\item[(iii)] $|Q_m(z)|\le C_1e^{-c_1m |z|^2}, \quad z\in K$.
\end{description}
\vskip-60pt
\be\label{descr}\ee
\vskip35pt

\noindent For some small
$\e>0$ consider
$R_n(z)Q_{\e n}(z)$. This is a rational function with poles in
$Z$ and at $\i$, and it will be important that the pole
at infinity coming from $Q_{\e n}$ is of order at most $\e n$.
We estimate this product on $\g_\tau$ as follows. Let  $z\in \g_\tau$
and let $0<\eta<1$ be given.
 If $|z|\le \sqrt{2\log (C_1)/c_1\e n}$,
then (\ref{ee}) and (ii)  yield
\[|R_n(z)Q_{\e n}(z)|\le \exp\Bigl(M_3\sqrt\tau 2\log (C_1)/c_1\e\Bigr)\|R_n\|_{\G},\]
and the right hand side is smaller than $(1+\eta)\|R_n\|_{\G}$ if
$\tau<(\eta c_1\e/4M_3\log C_1)^2$.
On the other hand, if
$|z|> \sqrt{2\log (C_1)/c_1\e n}$,
then (\ref{ee}) and (iii) give
\be |R_n(z)Q_{\e n}(z)|\le
\|R_n\|_{\G} C_1\exp\Bigl(n M_3\sqrt \tau|z|^2
-c_1\e n |z|^2\Bigr).\label{ui}\ee
For $\sqrt \tau<c_1 \e/2M_3$ the exponent is at most
\[-n (c_1/2)\e |z|^2\le \log (1/C_1)\]
so  in this case we have
\be |R_n(z)Q_{\e n}(z)|\le \|R_n\|_{\G}.\label{kl}\ee

What we have shown is that
\be \|R_nQ_{\e n}\|_{\gamma_\tau}\le (1+\eta)\|R_n\|_{\G}\label{hah}\ee
if $\tau$ is small, say $\tau<\tau_\eta^*$. Fix such a $\tau$. The corresponding
$\g_\tau$ is an analytic arc, so we can apply (\ref{mainarc}) to it and
to the rational function $R_nQ_{\e n}$ (recall that (\ref{mainarc}) has already been
proven for analytic arcs in Section \ref{sectanalarcs}). It follows that
\be |(R_nQ_{\e n})'(0)|\le (1+o(1)) \|R_nQ_{\e n}\|_{\g_\tau}
\max\left(\sum_a{}'\frac{\partial g_{\C\setm \g_\tau}(0,a)}{\partial {\bf n}_+},
\sum_a{}'\frac{\partial g_{\C\setm \g_\tau}(0,a)}{\partial {\bf n}_-}\right),
\label{kl1}\ee
where now $\sum_a'$ means that the summation is for the poles of $R_nQ_{\e n}$,
i.e. for the poles of $R_n$ as well as for the at most $\e n$ poles $a=\i$
that possibly come from $Q_{\e n}$. Note that some of the poles may be cancelled
in $R_nQ_{\e n}$, but the inequality
\be \sum_a{}'\frac{\partial g_{\C\setm \g_\tau}(0,a)}{\partial {\bf n}_\pm}
\le \sum_a\frac{\partial g_{\C\setm \g_\tau}(0,a)}{\partial {\bf n}_\pm}
+\e n\frac{\partial g_{\C\setm \g_\tau}(0,\i)}{\partial {\bf n}_\pm},\label{dfg}\ee
(where on the right the summation is only on the original poles of $R_n$)
holds in that case, as well.
For the first sum on the right we use (\ref{normal11})  and
for the second sum Lemma \ref{poscor} to conclude
\be \sum_a{}'\frac{\partial g_{\C\setm \g_\tau}(0,a)}{\partial {\bf n}_\pm}
\le (1+\eta)\sum_a\frac{\partial g_{\C\setm \G}(0,a)}{\partial {\bf n}_\pm}
+C_2\e n\label{ET}\ee
with some $C_2$ that depends only on $\G$. Since the
sum on the right of (\ref{ET}) is $\ge c_2n$ with some fixed $c_2>0$
again by Lemma \ref{poscor}, we obtain
from (\ref{hah}) and (\ref{kl1})
\bean  |(R_nQ_{\e n})'(0)|&\le& (1+o(1)) (1+\eta)^2\|R_n\|_{\G}
(1+C_2\e/c_2)\times \\[10pt]
&&\max\left(\sum_a\frac{\partial g_{\C\setm \G}(0,a)}{\partial {\bf n}_+},
\sum_a\frac{\partial g_{\C\setm \G}(0,a)}{\partial {\bf n}_-}\right).\eean

In view of $Q_{\e n}(0)=1$, on the left
\[(R_nQ_{\e n})'(0)=R_n'(0)+R_n(0)Q_{\e n}'(0),\]
and for the second term we get again from (\ref{mainarc}) (known
for the analytic arc $\g_\tau$ by Section \ref{sectanalarcs})
 and from
 $\|Q_{\e n}\|_{\g_\tau}\le 1$
\[|R_n(0)Q_{\e n}'(0)|\le (1+o(1)) \|R_n\|_{\G}
n\e\max\left(\frac{\partial g_{\C\setm \g_\tau}(0,\i)}{\partial {\bf n}_+},
\frac{\partial g_{\C\setm \g_\tau}(0,\i)}{\partial {\bf n}_-}\right),\]
and we can again apply (\ref{normal}) to the right hand side.
If we use again Lemma \ref{poscor} as before, we finally obtain
\[|R_n'(0)|\le (1+o(1)) (1+\eta)^2(1+C_3\e)\|R_n\|_{\G}
\max\left(\sum_a\frac{\partial g_{\C\setm \G}(0,a)}{\partial {\bf n}_+},
\sum_a\frac{\partial g_{\C\setm \G}(0,a)}{\partial {\bf n}_-}\right)\]
with some constant $C_3$ independent of $\e$ and $\eta$.
Now this is true for all $\e,\eta>0$
so  the claim (\ref{mainarc}) follows.

We shall not prove the last statement concerning the uniformity of the
estimate, for the argument is very similar to the one
given in the proof of \cite[Theorem 1]{TotikMarkov}.\endproof
\bigskip

\sect{Proof of Theorem \ref{thcurve}}
In this section we prove the inequality
(\ref{maincurve})  for $C^2$ curves. Recall
that in Section \ref{sectanalcurve} the inequality  (\ref{maincurve})  has already been proven for
analytic curves, which was the basis of all
subsequent results. In the present section we show
how (\ref{maincurve})  for $C^2$ curves can be deduced from the inequality
(\ref{mainarc})
for $C^2$ arcs.

Thus, let $\G$ be a positively oriented $C^2$ smooth Jordan curve
and $z_0$ a point on $\G$. Let $w_0\not =z_0$ be another point of
$\G$ (think of $w_0$ as lying ``far" from $z_0$), and for
$m=1,2,\ldots$ let $w_m\in \G$ be the point on $\G$ such that
the arc $\ov{w_0w_m}$ (in the orientation of $\G$) is of length $1/m$.
Such a $w_m$ exists and the arc $\overline{w_0w_m}$  does not contain
$z_0$ for all sufficiently large $m$, say for $m\ge m_0$.
Remove now the (open) arc $\overline{w_0w_m}$ from $\G$ to get
the Jordan arc $\G_m=\G\setm \overline{w_0w_m}$. We can apply (\ref{mainarc}) to
this $\G_m$, and what we are going to show is that
the so obtained inequality proves (\ref{maincurve})
as $m\to\i$.

Let $a\in G_-\cap Z$. We show that, as $m\to\i$,
\be \frac{\partial g_{\C\setm \G_m}(z_0,a)}{\partial {\bf n}_-}
\to \frac{\partial g_{G_-}(z_0,a)}{\partial {\bf n}_-}
\label{b0}\ee
and
\be \frac{\partial g_{\C\setm \G_m}(z_0,a)}{\partial {\bf n}_+}
\to 0,
\label{b1}\ee
uniformly in $a\in G_-\cap Z$.
Indeed, since $\G_m\subset\G_{m+1}$, the Green's functions $g_{\C\setm \G_m}(z,a)$
decrease as $m$ increases. Furthermore,
$g_{\C\setm \G_{m_0}}(z,a)$ is continuous at $w_0$, so
for every $\e>0$ there is an $m_\e$ such that  for $z\in \overline{w_0w_{m_\e}}$ we have
$g_{\C\setm \G_{m_0}}(z,a)<\e$. In view of Corollary \ref{corgreencomp} this $m_\e$ can
be the same for all $a\in Z\cap G_-$ since the Green's functions
$g_{\C\setm \G_{m_0}}(z,a)$ with respect to different $a\in Z\cap G_-$
are comparable inside a Jordan curve $\s$ that encloses $\G_{m_0}$.
This then implies for $m\ge m_\e$
and $z\in \overline{w_0w_m}$
\be 0<g_{\C\setm \G_{m}}(z,a)\le g_{\C\setm \G_{m_0}}(z,a)<\e.\label{b2}\ee
Thus, for $m\ge m_\e$ the function $g_{\C\setm \G_{m}}(z,a)-g_{G_-}(z,a)$
is positive and harmonic in $G_-$, and on the boundary of $G_-$ it is either
0 or $<\e$, so by the maximum principle it is $<\e$ everywhere in
the closure $\overline {G_-}$. Let now $a_0\in G_+$ be fixed, i.e.
$a_0$ lies in the outer domain to $\G$, and let $I\subset \G$ be
a subarc of $\G$ which does not contain $z_0$ and which contains
$\overline{w_0w_{m_0}}$ in its interior, and set $\G_I=\G\setm I$. Then
$g_{\C\setm \G_I}(z,a_0)$ has a strictly positive lower bound $c_0$
on $\overline{w_0w_{m_0}}$ (note that this arc lies inside
the domain $\C\setm \G_I$), therefore, in view of (\ref{b2}),
we have
\be  0<g_{\C\setm \G_{m}}(z,a)-g_{G_-}(z,a)<\frac{\e}{c_0}g_{\C\setm \G_I}(z,a_0)
\label{b3}\ee
on the boundary of $G_-$ provided $m\ge m_\e$. By the maximum principle
this inequality then holds throughout $G_-$ (note that both sides are
harmonic there), and hence we have for $m\ge m_\e$
\be  0\le \frac{\partial g_{\C\setm \G_{m}}(z_0,a)}{\partial {\bf n}_-}
-\frac{\partial g_{G_-}(z_0,a)}{\partial {\bf n}_-}
\le \frac{\e}{c_0}\frac{\partial g_{\C\setm \G_I}(z_0,a_0)}{\partial {\bf n}_-},\label{cec}\ee
and upon letting $\e\to 0$ we obtain (\ref{b0}).

The proof of (\ref{b1}) is much the same, just work now in the
exterior domain $G_+$, and use the reference Green's function
$g_{\C\setm \G_I}(z,b_0)$ with $b_0$ lying in the bounded
domain $G_-$. In this case $g_{\C\setm \G_{m}}(z,a)$ is harmonic
in $G_+$ for $a\in G_-\cap Z$ and (\ref{b3}) takes
the form
\[  0<g_{\C\setm \G_{m}}(z,a)<\frac{\e}{c_0}g_{\C\setm \G_I}(z,b_0),\]
from where the conclusion (\ref{b1}) can be made
as before.

For poles $a$ lying outside $\G$ we have similarly
\be \frac{\partial g_{\C\setm \G_m}(z_0,a)}{\partial {\bf n}_+}
\to \frac{\partial g_{G_+}(z_0,a)}{\partial {\bf n}_+}
\label{b00}\ee
and
\be \frac{\partial g_{\C\setm \G_m}(z_0,a)}{\partial {\bf n}_-}
\to 0,
\label{b10}\ee
uniformly in $a\in G_+\cap Z$ as $m\to\i$.

After these preparations we turn to the proof of (\ref{maincurve}). Choose, for a large $m$, the Jordan arc
$\G_m$, and apply (\ref{mainarc}) to this Jordan arc and to the rational function $R_n$ in
Theorem \ref{thcurve}. Since $\|R_n\|_{\G_m}\le \|R_n\|_{\G}$, we obtain
\be |R_n'(z_0)|\le (1+o(1)) \|R_n\|_{\G}
\max\left(\sum_{a\in Z}\frac{\partial g_{\C\setm \G_m}(z_0,a)}{\partial {\bf n}_+},
\sum_{a\in Z}\frac{\partial g_{\C\setm \G_m}(z_0,a)}{\partial {\bf n}_-}\right)
\label{b54}\ee
where the $o(1)$ term may depend on $m$. In view of
(\ref{b0})--(\ref{b1}) and (\ref{b00})--(\ref{b10}) (use also (\ref{b41}) and (\ref{b410}))
\[\sum_{a\in Z}\frac{\partial g_{\C\setm \G_m}(z_0,a)}{\partial {\bf n}_+}
\le (1+o_m(1))\sum_{a\in Z\cap G_+}\frac{\partial g_{G_+}(z_0,a)}{\partial {\bf n}_+}
+o_m(1) n\]
and
\[\sum_{a\in Z}\frac{\partial g_{\C\setm \G_m}(z_0,a)}{\partial {\bf n}_-}
\le (1+o_m(1))\sum_{a\in Z\cap G_-}\frac{\partial g_{G_-}(z_0,a)}{\partial {\bf n}_-}
+o_m(1) n,\]
where $o_m(1)$ denotes a quantity that tends to 0 as $m\to\i$. These imply
that the maximum on the right of (\ref{b54}) is at most
\[(1+o_m(1))\max\left(\sum_{a\in Z\cap G_+}\frac{\partial g_{G_+}(z_0,a)}{\partial {\bf n}_+}
+o_m(1) n,\sum_{a\in Z\cap G_-}\frac{\partial g_{G_-}(z_0,a)}{\partial {\bf n}_-}
+o_m(1) n\right),\]
which is
\[(1+o_m(1))\max\left(\sum_{a\in Z\cap G_+}\frac{\partial g_{G_+}(z_0,a)}{\partial {\bf n}_+},
\sum_{a\in Z\cap G_-}\frac{\partial g_{G_-}(z_0,a)}{\partial {\bf n}_-}
\right)\]
because of (\ref{b41})--(\ref{b410}). Therefore, we obtain (\ref{maincurve}) from (\ref{b54}) by
letting $n\to\i$ and at the same time $m\to\i$ very slowly.

A routine check shows that the proof runs uniformly in  $z_0\in \G$
lying on any proper arc $J$ of $\G$. In fact, the proof gives that uniformity
provided the normal derivative on the right of (\ref{cec}) lies in between
two positive constants independently of $z_0\in J$, which can be easily proven using
the method of Proposition \ref{propgreencont} (which was  based on  the Kellogg-Warschawski theorem
and that  is uniform in $z_0$ in the given range). From here the uniformity of
(\ref{maincurve}) in $z_0\in \G$ follows by considering two such arcs  $J$
that together  cover $\G$. \endproof

\sect{Proof of (\ref{mainarck})}
In the proof of Theorem \ref{thMarkovk} we shall need
 (\ref{mainarck}) which we verify in this section.
 The proof uses  induction on $k$, the $k=1$ case is covered by  Theorem  \ref{tharc}.

Let $R_n$ and $J$ as in  (\ref{mainarck}).
First of all we remark that by \cite[Theorem 7.1]{TotikTrans},
$g_{\C\setm \G}(z,\i)$ is  H\"older 1/2 continuous: for all $z\in \C$
\[ g_{\C\setm \G}(z,\i)\le M{\rm dist}(z,\G)^{1/2}\]
with some constant $M$. This combined with Corollary \ref{corgreencomp} (just apply it
to $\g_0=\G$) shows that all $g_{\C\setm \G}(z,a)$, $a\in Z$, are
uniformly  H\"older 1/2 equicontinuous:
\[ g_{\C\setm \G}(z,a)\le M_1{\rm dist}(z,\G)^{1/2}, \qquad a\in Z,\  {\rm dist}(z,\G)\le d,\]
with some constants $M_1$ and $d>0$. If we use also (\ref{aa}), then we obtain
\[|R_n(z)|\le \|R_n\|_\G \exp\left(nM_1{\rm dist}(z,\G)^{1/2}\right), \qquad {\rm dist}(z,\G)\le d.\]
In particular, if $z_0\in \G$ and $C_{1/n^2}(z_0)$ is the circle
about $z_0$ of radius $1/n^2$, then for all $z\in C_{1/n^2}(z_0)$ we have
$|R_n(z)|\le \|R_n\|_\G \exp(M_1)$. Thus, Cauchy's integral formula for the $k$-th derivative
at $z_0$ (written as a contour integral over $C_{1/n^2}(z_0)$) gives for large $n$
\[|R_n^{(k)}(z_0)|\le k!n^{2k}e^{M_1}\|R_n\|_\G,\]
and since this is true uniformly for all $z_0\in \G$,
\be \|R_n^{(k)}\|_\G\le C_kn^{2k}\|R_n\|_\G\label{pom}\ee
follows with some $C_k$.

Let
\[V(u )=
\max\left(\sum_{i=0}^mn_i\frac{\partial g_{\C\setm \G}(u,a_i)}{\partial {\bf n}_+},
\sum_{i=0}^mn_i\frac{\partial g_{\C\setm \G}(u,a_i)}{\partial {\bf n}_-}\right).\]
We shall need the following equicontinuity property of these $V(u)$:
\be V(v)\le (1+\e)V(z_0) \qquad \mbox{if $z_0\in J$ and $|v-z_0|<\d,\ v\in \G$},\label{eqi}\ee
with some $\e$ that tends to 0 as $\d\to0$. It is clear that this follows if
we prove the continuity for each term in $V(u)$, for example, if we show that
\be \frac{\partial g_{\C\setm \G}(v,a)}{\partial {\bf n}_-}\le
(1+\e)\frac{\partial g_{\C\setm \G}(z_0,a)}{\partial {\bf n}_-}\label{zao}\ee
if $z_0\in J$ and $|v-z_0|<\d$ where  $\e$ tends to 0 as $\d\to0$.
If $\f$ is a conformal map from the unit disk onto $\C\setm \G$ that maps
0 into $a$, then, just
as in (\ref{37}), we have
\be \frac{\partial g_{\C\setm \G}(v,a)}{\partial {\bf n}_-}=\frac{1}{|\f'(\f^{-1}(v)|},\label{370}\ee
with the understanding that of the two pre-images $\f^{-1}(v)$ of $v$, in this
formula we select the one that is mapped to the left side of $\G$ by $\f$.
A relatively simple localization (just open up the arc $\G$ to a $C^2$ Jordan curve as in
Section \ref{sectanalarcs})
of the  Kellogg-Warschawski theorem (\cite[Theorem 3.6]{Pommerenke})
 shows that $\f'$ is positive and continuous
away from the pre-images of the endpoints of $\G$. This  implies (\ref{zao}) in view
of (\ref{370}).

Suppose now that the claim in (\ref{mainarck})
is true for a $k$ and for all
subarcs $J\subset \G$ that do not contain either of the endpoints of $\G$. For such a subarc select
a subarc $J\subset J^*$ such that $J^*$ has no common endpoint either with
$J$ or with $\G$.
For a $z_0\in J$ let $Q(v)=Q_{n^{1/3},z_0}(v)$ be as in (i)--(iii) of (\ref{descr})
with 0 replaced by $z_0$ and $K$ replaced by $\G$.
So this is a polynomial of degree at most $ n^{1/3}$ such that $Q(z_0)=1$, $\|Q\|_\G\le 1$ and
if $v\in \G$, then
\be |Q(v)|\le C_1e^{-c_1n^{1/3} |v-z_0|^{2}}.\label{21}\ee
Because of the uniform $C^2$ property of $\G$ relatively
simple consideration shows that here  the constants $C_1,c_1$ are independent of $z_0\in J$.

Consider any $\d>0$ such that the intersection of $\G$ with the $\d$-neighborhood of
$J$ is part of $J^*$, and set
$f_{k,n,z_0}(v)=R_n^{(k)}(v)Q(v)$. On $\G$ for this we have the bound
\[O(n^{2k})\exp(-c_1n^{1/3} \d^{2})\|R_n\|_\G=o(1)\|R_n\|_\G\]
outside the $\d$-neighborhood of $z_0 $ (see (\ref{pom}) and (\ref{21})).
In the $\d$-neighborhood of any $z_0\in J$ we have,
by $\|Q\|_\G\le 1$ and by the induction hypothesis applied to $R_n$ and to the arc $J^*$,
\bean |f_{k,n,z_0}(v)|&\le& (1+o(1))\|R_n\|_\G V(v)^k\\
&\le& (1+o(1))(1+\e)^k\|R_n\|_\G V(z_0 )^k,\eean
where $\e\to0$ as $\d\to 0$ in view of (\ref{eqi}).
Therefore, $f_{k,n,z_0}(v)$ is a rational function in $v$
of total degree at most $n+n^{1/3}+mk$ (see below) for which
\[\|f_{k,n,z_0}\|_\G\le (1+o(1))\|R_n\|_\G V(z_0 )^k,\]
where $o(1)\to 0$ uniformly as $n\to\i$. The poles of $f_{k,n,z_0}$ agree with the poles $a_i$
of $R_n$ with a slight modification: for $a_i\not =\i$ the order of $a_i$
in $f_{k,n,z_0}$ is at most $n_i+k$ (see the form (\ref{spec}) of $R_n$), while for $a_0=\i$
the order of $a_0$ is at most $n_0-k$ plus  at most $n^{1/3}$
coming from $Q$. Upon applying Theorem \ref{tharc} to the rational function
$f_{k,n,z_0}$ we obtain (see also (\ref{b41}) and (\ref{b410}))
\bean  |f_{k,n,z_0}'(z_0 )|&\le& (1+o(1))\|R_n\|_\G V(z_0 )^k \times\\
&&\left(V(z_0)+O(mk)+n^{1/3}
\max\left(\frac{\partial g_{\C\setm \G}(z_0,\i)}{\partial {\bf n}_+},
\frac{\partial g_{\C\setm \G}(z_0,\i)}{\partial {\bf n}_-}\right)\right).\eean
In view of (\ref{b41})--(\ref{b410}) $V(z_0)$ is much larger (of size $n$)
than the last two
terms on the right (which are together of size $O(n^{1/3})$ if $z_0$ stays away
from the endpoints of $\G$), hence it follows that
\be |f_{k,n,z_0}'(z_0 )|\le (1+o(1))\|R_n\|_\G V(z_0 )^{k+1} .\label{fhg}\ee

Since (recall that $Q(z_0)=1$)
\[f_{k,n,z_0}'(z_0 )=R_n^{(k+1)}(z_0 )+R_n^{(k)}(z_0 )Q'(z_0),\]
and the second term on the right is
$O(n^{2/3})O(n^k)\|R_n\|_\G$ by
the induction assumption and by (\ref{pom}) applied to $Q$ with $k=1$ rather than
to $R_n$,  we
can conclude (\ref{mainarck}) for $k+1$ from (\ref{fhg}).

From how we have derived this, it follows that this estimate is uniform in $z_0\in J$.\endproof

\sect{The Markov-type inequality for higher derivatives}
In this section we prove the first part of Theorem \ref{thMarkovk} (the
sharpness will be handled in Section \ref{sectsharp}).
The proof uses the symmetrization
technique of \cite{TotikMarkov}.
It is sufficient to prove (\ref{mainMarkovaltk}).

First of all we remark that the limits defining $\O_a(A)$ in (\ref{Odef}) exist and are equal
for the choices ${\bf n}_\pm$.
Indeed, let $\f_a(z)=1/(z-a)$ be the fractional linear transformation considered
before. Then
\[g_{\C\setm \G}(z,a)=g_{\C\setm \f_a(\G)}(\f_a(z),\i),\]
so for a $z\in \G$ we have
\[\frac{\partial g_{\C\setm \G}(z,a)}{\partial {\bf n}_\pm}=
\frac{\partial g_{\C\setm \f_a(\G)}(\f_a(z),\i)}{\partial {\bf n}_\pm}|\f_a'(z)|,\]
and it has been verified in the proof of  \cite[Theorem 2]{TotikMarkov} that, as $w\to\f_a(A)$, $w\in \f_a(\G)$,
\[\sqrt{|w-\f_a(A)|}\frac{\partial g_{\C\setm \f_a(\G)}(w,\i)}{\partial {\bf n}_\pm}\]
have equal limits, call them $\O_\i(\f_a(\G),\f_a(A))$,
for both  choices of $+$ or $-$.
Since, as $z\to A$, $z\in \G$, we have
$|\f_a(z)-\f_a(A)|=(1+o(1))|z-A||\f_a'(A)|$, it follows that,
indeed,
the limits
\[\lim_{z\to A,\ z\in \G} \sqrt{|z-A|}\frac{\partial g_{\C\setm \G}(z,a)}{\partial {\bf n}_\pm}=
\O_\i(\f_a(\G),\f_a(A))\sqrt{|\f_a'(A)|}\]
exist and are the same for the $+$ or $-$ choices.

Next, we prove the required inequality at the endpoint $A$.
We may assume that $A=0$.
Let
\[\G^*=\{z\sep z^2\in \G\}.\]
This is a Jordan arc symmetric with respect to the origin.
It is not difficult to prove (see \cite[Appendix 2]{TotikMarkov})
that $\G^*$ has $C^2$ smoothness.

Let $R_n$ be a rational function of degree at most
$n$ of the form (\ref{spec}), and set ${\cal R}_{2n}(z)=R_n(z^2)$. This is a rational function
which has $2n$ poles $\pm \sqrt {a_i}$, where $a_i$ runs through
the poles of $R_n$ (here $\pm \sqrt {a_i}$ denote the two possible
values of $\sqrt {a_i}$ with the understanding
that if $a_0=\i$, then both values $\pm\sqrt{a_0}$
is $\i$).
If we apply (\ref{mainarck}) to $\G^*$ and to the rational function ${\cal R}_{2n}$, then
we get
\be |{\cal R}_{2n}^{(2k)}(0)|\le (1+o(1))M^{2k}\|{\cal R}_{2n}\|_{\G^*},\label{aa00}\ee
where
\be M=\max_\pm \sum_{i=0}^m n_i\left\{\frac{\partial g_{\C\setm \G^*}(0,\sqrt {a_i})}{\partial {\bf n}_\pm}
+\frac{\partial g_{\C\setm \G^*}(0,-\sqrt {a_i})}{\partial {\bf n}_\pm}\right\}.
\label{b36}\ee
For $a\not=\i$
\[g_{\C\setm \G}(z^2,a)=g_{\C\setm \G^*}(z,\sqrt a)+g_{\C\setm \G^*}(z,-\sqrt a),\]
hence for $z\not=0$ we have
\be \frac{\partial g_{\C\setm \G^*}(z,\sqrt a)}{\partial {\bf n}_\pm(z)}
+\frac{\partial g_{\C\setm \G^*}(z,-\sqrt a)}{\partial {\bf n}_\pm(z)}=
\frac{\partial g_{\C\setm \G}(z^2,a)}{\partial {\bf n}_\pm(z^2)}|2z|
\label{vb}\ee
(with possibly replacing ${\bf n}_\pm$ by ${\bf n}_\mp$ on the right),
which implies
\be \frac{\partial g_{\C\setm \G^*}(0,\sqrt a)}{\partial {\bf n}_\pm}
+\frac{\partial g_{\C\setm \G^*}(0,-\sqrt a)}{\partial {\bf n}_\pm}
=2\lim_{w\to 0}\frac{\partial g_{\C\setm \G}(w,a)}{\partial {\bf n}_\pm(w)}\sqrt{|w|}
=2\O_a(A).
\label{hhh}\ee
For $a=\i$ the corresponding calculation is
\[g_{\C\setm \G^*}(z,\i)=\frac12g_{\C\setm \G}(z^2,\i),\qquad \frac{\partial g_{\C\setm \G^*}(z,\i)}{\partial {\bf n}_\pm(z)}=
\frac12\frac{\partial g_{\C\setm \G}(z^2,\i)}{\partial {\bf n}_\pm(z^2)}|2z|,\]
and so
\be \frac{\partial g_{\C\setm \G^*}(0,\i)}{\partial {\bf n}_\pm}
=\lim_{w\to 0}\frac{\partial g_{\C\setm \G}(w,\i)}{\partial {\bf n}_\pm(w)}\sqrt{|w|}
=\O_\i(A).
\label{hhh1}\ee
Thus, the $M$ in (\ref{b36}) is exactly
\be 2\left(\sum_{i=0}^mn_i\O_{a_i}(A)\right).\label{fffh}\ee

In what follows we shall also need that the quantities $\O_{a_i}(A)$ are finite and positive,
which is immediate from (\ref{hhh}) and Lemma \ref{poscor} (this latter applied
to $\g_\tau=\g_0=\G$).

Now we use Fa\`{a} di Bruno's formula \cite{FB} (cf. \cite[Theorem 1.3.2]{Krantz})
\be  (S(F(z)))^{(2k)}=
\sum_{\nu_j} \frac{(2k)!}{\prod_{j=1}^{2k} \nu_j!(j!)^{\nu_j}}
S^{(\nu_1+\cdots+\nu_{2k})}(F(z))\prod_{j=1}^{2k} \left(F^{(j)}(z)\right)^{\nu_j},
\label{fd}\ee
where the summation is for all nonnegative integers $\nu_1,\ldots,\nu_{2k}$
for which $\nu_1+2\nu_2+3\nu_3+\cdots+2k\nu_{2k}=2k$, and where $0^0$ is defined
to be 1 if it occurs on the right. Apply this with $S=R_n$
and $F(z)=z^2$ at $z=0$:
\bean   {\cal R}_{2n}^{(2k)}(0)&=&{(R_{n}(F(z)))^{(2k)}\restr{z=0}}\\
&=&
\sum_{\nu_j} \frac{(2k)!}{\prod_{j=1}^{2k} \nu_j!(j!)^{\nu_j}}
R_{n}^{(\nu_1+\cdots+\nu_{2k})}(0)\prod_{j=1}^{2k} \left(F^{(j)}(0)\right)^{\nu_j}\\
&=&\frac{(2k)!}{k!2^k}R_{n}^{(k)}(0)2^k\eean
(use that $F^{(j)}(0)=0$ unless $j=2$ and then $F^{(2)}(0)=2$).
Hence, in view of (\ref{aa00}) we obtain
\be |R_{n}^{(k)}(0)|\le (1+o(1))\frac{2^k}{(2k-1)!!}
\left(\sum_{i=0}^m n_i\O_{a_i}(A)\right)^{2k}\|R_{n}\|_{\G},
\label{io1}\ee
where we also used that $\|{\cal R}_{2n}\|_{\G^*}=\|R_{n}\|_{\G}$.
This proves gives the correct bound for the $k$-th derivative
at the endpoint $A$.

\bigskip

So far we have verified (\ref{io1}), which is the claim (\ref{mainMarkovaltk}),
but only at the endpoint $A=0$ of the arc $\G$. We can reduce
the Markov type inequality (\ref{mainMarkovaltk}) to this special case.
To achieve that let us denote $\O_a(A)$ for the
arc $\G$ by $\O_a(\G,A)$.
If $z\in \G$ is close to $A$, then consider the subarc $\G_z$
which is the arc of $\G$ from $z$ to $B$ (recall that $B$ is the other
endpoint of $\G$ different from $A$), so the endpoints of $\G_z$
are $B$ and $z$. It is easy to see that the preceding proof
of (\ref{io1}) was uniform in the sense that it holds uniformly
for all $\G_z$, $z\in \G$, $|z-A|\le |B-A|/2$
(see the proofs of Theorem 3 and Appendix 1 in \cite{TotikMarkov}), therefore we obtain (replace in (\ref{io1})
$A$ by $z$)
\be |R_n^{(k)}(z)|\le (1+o(1))\frac{2^k}{(2k-1)!!}
\left(\sum_{i=0}^mn_i\O_{a_i}(\G_z,z)\right)^{2k}\|R_{n}\|_{\G_z},
\label{oo1}\ee
where
now the quantity $\O_a(\G_z,z)$ must be taken with respect to $\G_z$, rather than
with respect to $\G$. Since on the right
\[ \|R_n\|_{\G_z}\le  \|R_n\|_{\G},\]
all what remains to prove is that
\be \lim_{z\to A,\ z\in \G}\O_{a_i}(\G_z,z)\to \O_{a_i}(\G,A)\label{ll}\ee
for each $a_i$, $i=0,1,\ldots,m$, as $z\to A$.
Indeed, then we obtain from (\ref{oo1}) and from the fact that,
as has been mentioned before, the $\O_{a_i}(A)$ quantities
 are finite and positive, that
for any $\e>0$
\be  |R_n^{(k)}(z)|\le (1+\e) \frac{2^k}{(2k-1)!!} \left(\sum_{i=0}^mn_i\O_{a_i}(\G_z,z)\right)^{2k}\|R_n\|_{\G},\label{oo4}\ee
if $z\in \G$ lies sufficiently close to $A$, say $|z-A|\le \d$,
and $n$ is sufficiently large.
On the other hand, (\ref{mainarck}) shows that
$R_n^{(k)}(z)=O(n^k)$ on subsets of $\G$ lying away from the endpoints $A,B$,
in particular this is true
for $z\in U$, $|z-A|\ge \d$.
Now this and (\ref{oo4}) prove the theorem. So it is enough to prove
(\ref{ll}).

 (\ref{ll}) has been verified for $a_i=\i$ in the proof of
 \cite[Theorem 3]{TotikMarkov}. To get it for other $a_i$ just
 apply the mapping $\f_{a_i}(z)=1/(z-a_i)$ as before to reduce
 it to the $a_i=\i$ special case. The reader can easily fill in the
 details.\endproof

\sect{Proof of the sharpness}\label{sectsharp}
In this section we prove Theorems \ref{thcurvesh},
\ref{tharcsh} and the second part of Theorem \ref{thMarkovk}.

We shall first give the proof for Theorem \ref{thcurvesh}.
The proof of Theorem \ref{tharcsh} can be reduced to Theorem \ref{thcurvesh}
by attaching a suitable lemniscate as in the proof of Theorem \ref{thMarkovk},
so we skip it (actually,
a complete proof will be given as part of the proof in Section \ref{sectsh}
for rational functions
of the form (\ref{spec}) with fixed poles). However, the sharpness in
Theorem  \ref{thMarkovk} requires a different approach which will be given in
Section \ref{sectsh}.

\subsection{Proof of Theorem \ref{thcurvesh}}\label{sectthcurvesh}

The idea is as follows. On the unit circle,
we use some special rational functions (products of Blaschke factors)
for which  the Borwein-Erd\'elyi inequality (Proposition \ref{bepr})
is sharp. Then we transfer that back to $\Gamma$
 and approximate the transformed function with rational functions. In other words, we
reverse the reasoning in Section \ref{sectanalcurve} and
do the  ``reconstruction step'' in the  ``opposite direction''.

Recall that $\mathbf{D}=\{v\sep |v|<1\}$ and
$\mathbf{D}_+=\{v\sep |v|>1\}\cup\{\infty\}$, and
denote by $B(a,v)=\frac{1-\overline{a}v}{v-a}$ the (reciprocal) Blaschke factor with pole at $a$.

First, we state cases when we have equality in Proposition \ref{bepr}.
\begin{prop}
\label{prp:BEsharpness}
Suppose $h$ is a (reciprocal) Blaschke product with all poles either inside or outside the unit circle,
that is,
$h(v)=\prod_{j=1}^n B\left(\alpha_j,v\right)$ where all $\alpha_j\in\mathbf{D}$,
or $h(v)=\prod_{j=1}^n B\left(\beta_j,v\right)$ where all $\beta_j\in\mathbf{D}_+$.
Then
\[
\left|h'(1)\right|
=
\left\|h\right\|_{\unitcircle}
\max\left(
\sum_{\alpha_j}
\frac{\partial g_{\mathbf{D}}\left(1,\alpha_j\right)}{\partial \mathbf{n}_-},
\sum_{\beta_j}
\frac{\partial g_{\mathbf{D}_+}\left(1,\beta_j\right) }{\partial \mathbf{n}_+}
\right).
\]
\end{prop}
This proposition is contained in the  Borwein-Erd\'elyi theorem as stated in \cite{BE} pp. 324-326.

\medskip

First, we consider the case when
\begin{equation}
\Gamma \mbox{ is analytic and $Z\cap G_-\not=\emptyset$ } ,
\end{equation}
where, as always, $G_-$ is the interior domain determined by $\G$.

Fix $z_0\in \G$, and let, as in Section \ref{sectconf}, $\F_1$ be the conformal map from the unit disk
onto the interior domain $G_-$ such that $\F_1(1)=z_0$, $|\F_1'(1)|=1$.
As has been discussed there, this $\F_1$ can be extended to a disk $\{v\sep |v|<r_1\}$ with some
$r_1>1$.

Let $\alpha_1,\ldots,\alpha_n$ be $n$ (not necessarily different) points from
$\Phi_1^{-1}\left(ZÂ\cap G_-\right)$, and let
\[
h_n(v):=\prod_{j=1}^n
B\left(\alpha_j,v\right),
\]
for which $\left\|h_n\right\|_{\unitcircle}=1$.
Now we  ``transfer''  $h_n$ to $G_-$ by considering $h_n\left(\Phi_1^{-1}(z)\right)$.
If $f_{1,n}(z)$ is the sum of the principal parts of $h_n\left(\Phi_1^{-1}(z)\right)$
(with $f_{1,n}(\i)=0$), then
\[
\varphi_{e}(z):=h_n\left(\Phi_1^{-1}(z)\right) - f_{1,n}(z)
\]
is analytic in $G_{1}^+:=\left\{ \Phi_1(v)\sep |v|<r_{1}\right\}$.
Since $h_n$ is at most 1 in absolute value outside the unit disk,
it follows from
Proposition \ref{propgrig} as in Section \ref{sectanalcurve}
that the absolute value of $\f_\e$ is $\le C\log n$ on $G_1^+$.
By Proposition \ref{propappr} (applied to $K=\{\F_1(v)\sep |v|\le \sqrt {r_1}\}$ and to $\tau=\partial G_1^+$)
there are polynomials $f_{2,\sqrt n}$ of degree at most $\sqrt n$
such that $f_{2,\sqrt n}(z_0)=\f_e(z_0)$, $f_{2,\sqrt n}'(z_0)=\f_e'(z_0)$
and
\be \|\f_e-f_{2,\sqrt n}\|_K\le C(\log n) q^{\sqrt n}\label{asd}\ee
with some $C$ and $q<1$.
Therefore, if we set
\[
f_n(z):=f_{1,n}(z)+f_{2,\sqrt n},
\]
then this is a rational function with poles in $Z\cap G_-$ of total degree $n$
and with one pole at $\i$ of order  $\le \sqrt n=o(n)$.
For it \[
\left|f_n'\left(z_0\right)\right|
=|(h_n(\F_1^{-1}))'(z_0)|=
\left|h_n'(1)\right|
\]
since $\left|\Phi_1'(1)\right|=1$.
Furthermore $\left\|h_n\right\|_{\unitcircle}=1$ (recall that $\mathbb{T}$ is the unit circle), so we obtain from (\ref{asd})
\begin{multline*}
\left\|f_n\right\|_\Gamma
= \left\| f_{1,n} + f_{2,\sqrt n}\right\|_\Gamma
=\left\| f_{1,n} + \varphi_e + f_{2,\sqrt n} -\varphi_e \right\|_\Gamma
\\=
\left\| h_n(\Phi_1^{-1}) + f_{2,\sqrt n} -\varphi_e \right\|_\Gamma
=1+O\left( (\log n )q^{\sqrt n}\right)=1+o(1).
\end{multline*}

We use Proposition  \ref{prp:BEsharpness} for $h_n$, hence
\begin{multline*}
\left|f_n'\left(z_0\right)\right|
=
\left|h_n'(1)\right|
 =
\left\|h_n\right\|_{\unitcircle}
\sum_{\a_j}
\frac{\partial g_{\mathbf{D}}\left(1,\a_j\right)}{\partial \mathbf{n}_-}
 \ge
(1-o(1))\left\|f_n\right\|_\Gamma
\sum_{\a_j}
\frac{\partial g_{\mathbf{D}}\left(1,\a_j\right)}{\partial \mathbf{n}_-}.
\end{multline*}
Here, by Proposition \ref{prop:green_deriv_unitdisk},
\begin{multline*}
\sum_{\a_j}
\frac{\partial g_{\mathbf{D}}\left(1,\a_j\right)}{\partial \mathbf{n}_-}
 =
\sum_{\a_j}
\frac{\partial g_{G_-}\left(z_0,\Phi_1(\a_j)\right)}{\partial \mathbf{n}_-}
\\=
\max\left(
\sum_{\a_j}
\frac{\partial g_{G_-}\left(z_0,\Phi_1(\a_j)\right)}{\partial \mathbf{n}_-}
,
\sqrt n
\frac{\partial g_{G_+}\left(z_0,\i\right)}{\partial \mathbf{n}_+}
\right),
\end{multline*}
where, in the last step, we used that the first term in the max is $\ge cn$ with some $c>0$
(see (\ref{b41})),
so the last equality holds for large $n$.

Summarizing, we have proven that if $\Gamma$ is an analytic Jordan curve, $Z\subset \C\setm \G$ is a closed set,
 such that $Z\cap G_-\not=\emptyset$,
 then there exist rational functions $R_{n,-}$ with poles at any prescribed locations $a_{1,n},\ldots,a_{n,n}\in Z\cap G_-$
and with a pole at $\i$ of order $o(n)$ such that
\begin{equation}
\left|R_{n,-}'\left(z_0\right)\right|
\ge \left(1-o(1)\right)
\left\Vert R_{n,-}\right\Vert_\Gamma
\sum_{a_{j,n}}
\frac{\partial g_{G_-}\left(z_0,a_{j,n}\right)}{\partial \mathbf{n}_-},
\label{ineq:anshpouts}
\end{equation}
where  $o(1)$ depends on $\Gamma$ and $Z$ only.

Similarly, if $\Gamma$ is still an analytic Jordan curve and $Z\cap G_+\not=\emptyset$,
 then the same assertion holds
for some rational functions $R_{n,+}$ with prescribed poles at $a_{j,n}\in Z\cap G_+$
and with a pole  of order $\le \sqrt n$ at some given point $\z_0$ inside $\G$:
\begin{equation}
\left|R_{n,+}'\left(z_0\right)\right|
\ge \left(1-o(1)\right)
\left\Vert R_{n,+}\right\Vert_\Gamma
\sum_{a_{j,n}}
\frac{\partial g_{G_+}\left(z_0,a_{j,n}\right)}{\partial \mathbf{n}_+}.
\label{ineq:anshpouts1}
\end{equation}
This follows by applying a suitable inversion: fix $\zeta_0\in G_-$
and apply the mapping $w=1/(z-\zeta_0)$.
We omit the details.

Now for analytic $\G$ Theorem \ref{thcurvesh} can be easily proven.
For simplicity assume that the $a_{j,n}$ are different and finite
(the following argument needs only simple modification if this is not
the case).
Suppose, for example, that for a given $n=1,2,\ldots$
\begin{equation}
\label{cond:outerbigger}
\sum_{a_{j,n}Â\in Z\cap G_-}
\frac{\partial g_{G_-}\left(z_0,a_{j,n}\right) }{\partialÂ\mathbf{n}_-}
\ge
\sum_{a_{j,n} \in Z\cap G_+}
\frac{\partial g_{G_+}\left(z_0,a_{j,n}\right) }{\partial \mathbf{n}_+}.
\end{equation}
Consider the poles $a_{j,n}$ that are in $G_-$, and denote by $R_-(z)$ a
rational function whose existence is established above
for these poles (if the number of the $a_{j,n}$ that are in $G_-$ is $N$,
then in the previous notation this $R_-$ is $R_{N,-}$, so  the number of poles of
$R_-$ in $G_-$ is $N$, and $R_-$ also has a pole of order
at most $\sqrt N$ at $\i$).
Next, for any given $\e>0$ write
\[f_{n,+}(z):=\varepsilon_n \sum_{a_{j,n} \in Z\cap G_+} \frac{1}{z-a_{j,n}}\]
 where $\varepsilon_n>0$ is so small that $\left\Vert f_{n,+}\right\Vert_{\Gamma}\le
 \varepsilon\left\Vert R_-\right\Vert_{\Gamma}$
and
$\left|f_{n,+}'\left(z_0\right)\right|Â\le \varepsilon \left|R_-'\left(z_0\right)\right|$.
It is easy to see that then $R_n(z):=R_-(z)+f_{n,+}(z)$ has poles at the prescribed points $a_{1,n},\ldots,a_{n,n}$ plus one additional
pole of order $\le \sqrt n$ at $\i$. Furthermore, it satisfies
\[\left|R_n'\left(z_0\right)\right|\ge (1-\varepsilon)^2(1-o(1))\|R_n\|_\G\sum_{a_{j,n}\in G_-}
\frac{\partial g_{G_-}(z_0,a_{j,n})}{\partial\mathbf{n}_-},\]
and, by the assumption (\ref{cond:outerbigger}),
the sum on the right is the same as the maximum in (\ref{maincurvesh}).

If \eqref{cond:outerbigger} does not hold (i.e. the
reverse inequality is true), then use the analogous $R_{+}$ ($=R_{n-N,+}$) and add to it a small
multiple of the sum of $1/(z-a_{j,n})$ with $a_{j,n}\in Z\cap G_-$.

Since in these estimates $\e>0$ is arbitrary, Theorem \ref{thcurvesh} follows for analytic $\G$.
\bigskip

If $\G$ is not analytic, only $C^2$ smooth, then we can do the following. Suppose for example,
that for an $n$ (\ref{cond:outerbigger}) is true. For $\e>0$ choose
an analytic Jordan curve, say  a lemniscate $L$, close to $\G$
such that $L\cap \Gamma=\{z_0\}$,
$L\setminus \{z_0\}$ lies in the interior of $\G$,   and
\be (1-\varepsilon) \frac{\partial g_{G_-}\left(z_0,\beta\right)}{\partial \mathbf{n}_-}\le
\frac{\partial g_{\C\setminus L}\left(z_0,\beta\right)}{\partial \mathbf{n}_-}\label{using}\ee
for all $\beta \in G_-\cap Z$. (Here  we used the shorthand notation $g_{\C\setm L}(z,a)$
for both $g_{{\rm Int}(L)}(z,a)$ when $a$ is inside $L$ and for $g_{{\rm Ext}(L)}(z,a)$ when $a$
is outside $L$, where ${\rm Int}(L)$ and ${\rm Ext}(L)$ denote the interior and
exterior domains to $L$.)
The existence of $L$ follows from the sharp form of Hilbert's lemniscate theorem in
\cite[Theorem 1.2]{NagyTotik1}
when $\b=\i$. For other $\b$ use fractional linear transformations to move the pole $\b$ to $\i$,
see the formula (\ref{iiu}) below, as well as the reasoning there.

Now construct $R_n$ for this $L$ as before, and multiply it by a polynomial $Q=Q_{n^{7/8}}$ of degree
at most $n^{7/8}$ such that
 $Q(1)=1$, $\|Q\|_{\G}\le 1$,  and with some constants $c_0,C_0>0$
\[
|Q(z)|\le C_0\exp(-c_0 n^{7/8}|z-z_0|^{3/2}),\qquad z\in \G.
\]
Such a $Q$ exists by \cite[Theorem 4.1]{TotikTrans}, and we have to consider
$R_nQ$ rather than $R_n$ because the norm of $R_n$ on $\G$ can be much larger
than its norm on $L$, and $Q$ brings that  norm down, namely
$\|R_nQ\|_\G\le (1+o(1))\|R_n\|_L$. Indeed, this is an easy consequence of (\ref{aa}) and
Proposition \ref{propgreencont}
(both applied to $L$ rather than $\G$) and the properties
of $Q$. Finally,  since $R_n$ proves Theorem \ref{thcurvesh} on $L$, relatively simple
argument shows that $R_nQ$ verifies it on $\G$. The reader can easily fill in
the details.\endproof

\subsection{Sharpness of the Markov inequality}\label{sectsh}
First we consider a $C^2$ Jordan curve $\g$ and a point $z_0\in \g$ on it.
Let $\e>0$. By the sharp form of the Hilbert lemniscate theorem \cite[Theorem 1.2]{NagyTotik1}
there is a Jordan curve $\s$ such that
\begin{itemize}
\item $\s$ contains $\g$ in its interior
except for the point $z_0$, where the two curves touch each other,
\item $\s$ is a lemniscate, i.e. $\s=\{z\sep |T_N(z)|=1\}$ for some
polynomial $T_N$ of degree $N$, and
\item
\be \frac{\partial g_{\C\setm \s}(z_0,\i)}{\partial {\bf n}_+}
\ge (1-\e)\frac{\partial g_{\C\setm \g}(z_0,\i)}{\partial {\bf n}_+},\label{xc}\ee
where the Green's functions $g_{\C\setm \g}(z_0,\i)$ and
$g_{\C\setm \s}(z_0,\i)$ are taken with respect to the outer domains
of $\g$ and $\s$.
\end{itemize}
We may assume that $T_N(z_0)=1$ and $T_N'(z_0)>0$.
The Green's function of the outer domain of $\s$ is $\frac1N \log|T_N(z)|$,
and its normal derivative is
\[\frac{\partial g_{\C\setm \s}(z_0,\i)}{\partial {\bf n}_+}=\frac1N |T_N'(z_0)|
=\frac1N T_N'(z_0).\]

Consider now, for all large $n$, the polynomials $S_n(z)=T_N(z)^{[n/N]}$, where
$[n/N]$ denotes integral part. This is a polynomial of degree at most $n$,
its supremum norm on $\s$ is 1, and
\[S_n'(z_0)=\left[\frac nN\right]T_N(z_0)^{[n/N]-1}T_N'(z_0)=n\frac{\partial g_{\C\setm \s}(z_0,\i)}{\partial {\bf n}_+}+O(1).\]
In a similar fashion,
\bean S_n''(z_0)&=&\left[\frac nN\right]\left(\left[\frac nN\right]-1\right)T_N(z_0)^{[n/N]-2}(T_N'(z_0))^2
+\left[\frac nN\right]T_N(z_0)^{[n/N]-1}T_N''(z_0)\\
&=&n^2\left(\frac{\partial g_{\C\setm \s}(z_0,\i)}{\partial {\bf n}_+}
\right)^2+O(n).\eean
Proceeding similarly, it follows that for any $j=1,2,\ldots$
\[S_n^{(j)}(z_0)=n^j\left(\frac{\partial g_{\C\setm \s}(z_0,\i)}{\partial {\bf n}_+}
\right)^j +O(n^{j-1}).\]
Thus, in view of (\ref{xc}), we may write
\be S_n^{(j)}(z_0)\ge (1-\e)^jn^j\left(\frac{\partial g_{\C\setm \g}(z_0,\i)}{\partial {\bf n}_+}
\right)^j +O(n^{j-1}),\label{snj}\ee
where, and in what follows, we use the following convention: if $A$ is a complex number and
$B$ is a positive number, then we write $A\ge B+O(n^s)$ if $A=C+O(n^s)$, where $C$ is a real
number with $C\ge B$. Note also that $\|S_n\|_\g\le\|S_n\|_\s=1$ by the
maximum principle.
\bigskip

Next, we need an analogue of these for rational functions with pole at a point $a$ that lies outside
$\g$. Consider the fractional linear transformation $\f_a(z)=\xi/(z-a)$, where $\xi$ is selected
so that $|\xi|=1$ and $\f_a'(z_0)>0$. The image of $\g$ under this transformation
is $\f_a(\g)$, and $g_{\C\setm \g}(z,a)=g_{\C\setm \f_a(\g)}(\f_a(z),\i)$. This latter relation
implies that
\be\frac{\partial g_{\C\setm \g}(z_0,a)}{\partial {\bf n}_+(z_0)}=
\frac{\partial g_{\C\setm \f_a(\g)}(\f_a(z_0),\i)}{\partial {\bf n}_+(\f_a(z_0))}\f_a'(z_0).\label{iiu}\ee
Now let $S_n$ be the polynomial constructed before, but this time for the
curve $\f_a(\g)$ and for the point $\f_a(z_0)$, and set
$S_{n,a}(z)=S_n(\f_a(z))$. This is a rational function with a pole
of order at most $n$ at $a$. Its norm on $\g$ is at most 1, and, in view of (\ref{snj}) (applied to $\f_a(\g)$),
\[S_{n,a}'(z_0)=S_n'(\f_a(z_0))\f_a'(z_0)
\ge (1-\e)n\frac{\partial g_{\C\setm \f_a(\g)}(\f_a(z_0),\i)}{\partial {\bf n}_+(\f_a(z_0))}\f_a'(z_0)
+O(1),\]
which can be written in the form
\[S_{n,a}'(z_0)\ge (1-\e)n\frac{\partial g_{\C\setm \g}(z_0,a)}{\partial {\bf n}_+}
+O(1)\]
in view of (\ref{iiu}). For the second derivative we have
\[S_{n,a}''(z_0)=S_n''(\f_a(z_0))(\f_a'(z_0))^2+S_n'(\f_a(z_0))\f_a''(z_0),\]
hence
\[S_{n,a}''(z_0)\ge (1-\e)^2n^2\left(\frac{\partial g_{\C\setm \g}(z_0,a)}{\partial {\bf n}_+}\right)^2
+O(n),\]
in view of (\ref{snj}) (applied to $\f_a(\g)$ and to the point $\f_a(z_0)$) and  (\ref{iiu}).
Proceeding similarly we obtain for all $j=1,2,\ldots$
\be S_{n,a}^{(j)}(z_0)\ge (1-\e)^jn^j\left(\frac{\partial g_{\C\setm \g}(z_0,a)}{\partial {\bf n}_+}
\right)^j +O(n^{j-1}).\label{snj1}\ee
\bigskip

Now let there be given a fixed number of different poles $a_0,\ldots,a_m$  in the exterior of $\g$ and
associated orders $n_0,\ldots,n_m$, where $a_0=\i$ (if we do not want
the point $\i$ among the poles, just set $n_0=0$). For the total degree $n=n_0+\cdots+n_m$
consider the rational function
\[U_n(z)=\prod_{i=0}^m S_{n_i,a_i}(z),\]
where we set $S_{n_0,a_0}=S_{n_0}$, with the
polynomial $S_{n_0}$ constructed in the first part of the proof. This is a rational
function with poles at the $a_i$'s and the order of $a_i$ is at most $n_i$.
Since
\[U_n^{(k)}(z_0)=
\sum_{j_0+\cdots +j_m=k} \frac{k!}{j_0!\cdots j_m!}S_{n_0,a_0}^{(j_0)}(z_0)\cdots
S_{n_m,a_m}^{(j_m)}(z_0),\]
we obtain from (\ref{snj1}) that
\bean U_n^{(k)}(z_0)&\ge &
\sum_{j_0+\cdots j_m=k} \frac{k!}{j_0!\cdots j_m!}\prod_{i=0}^m
\left((1-\e)^{j_i}n_i^{j_i}\left(\frac{\partial g_{\C\setm \g}(z_0,a_i)}{\partial {\bf n}_+}
\right)^{j_i} +O(n_i^{j_i-1})\right)\\
&&=\sum_{j_0+\cdots j_m=k} \frac{k!}{j_0!\cdots j_m!}\prod_{i=0}^m
\left((1-\e)^{j_i}n_i^{j_i}\left(\frac{\partial g_{\C\setm \g}(z_0,a_i)}{\partial {\bf n}_+}
\right)^{j_i}\right) +O(n^{k-1}),
\eean
so by the multinomial theorem (see e.g. \cite[Theorem 1.3.1]{Krantz})
\[ U_n^{(k)}(z_0)\ge (1-\e)^k\left(\sum_{i=0}^m n_i \frac{\partial g_{\C\setm \g}(z_0,a_i)}{\partial {\bf n}_+}
\right)^k+O(n^{k-1}).\]
Hence,
\be  |U_n^{(k)}(z_0)|\ge (1-\e)^k(1-o(1))\left(\sum_{i=0}^m n_i \frac{\partial g_{\C\setm \g}(z_0,a_i)}{\partial {\bf n}_+}
\right)^k\label{yut}\ee
in view of (\ref{b410}).

\bigskip
After these preparations we can prove the last statement in Theorem \ref{thMarkovk}. Let $\G$ be
a $C^2$ smooth Jordan arc and $a_0,\ldots,a_m$ be finitely many fixed poles
outside $\G$ with associated orders $n_0,\ldots,n_m$. We agree that $a_0=\i$,
and if we do not want
the point $\i$ among the poles, just set $n_0=0$. We may assume that the endpoint
$A$ of $\G$ is at the origin, and consider, as before, the curve
$\G^*=\{z\sep z^2\in \G\}$. We also consider the poles $\pm \sqrt{a_i}$, $i=0,\ldots,m$,
with associated orders $n_i$ with the agreement that if $n_0\not=0$, i.e.
the point $\i$ is among our poles, then $\pm\sqrt{\i}=\i$.

It is easy to see that there is a $C^2$ Jordan curve $\g$ such that
\begin{itemize}
\item $\g$ contains $\G^*$ in
its interior except for the point 0, where $\g$ and $\G^*$ touch each other,
\item all $a_i$ are outside $\g$,
\item
\be \frac{\partial g_{\C\setm \g}(0,a_i)}{\partial {\bf n}_+}\ge (1-\e)
\frac{\partial g_{\C\setm \G^*}(0,a_i)}{\partial {\bf n}_+}\quad \mbox{for all $i$}.\label{sdf}\ee
\end{itemize}
Indeed, all we need to do is to select $\g$ sufficiently close to
$\G^*$ and to have at 0 curvature close to the curvature of $\G^*$,
see e.g. \cite{NagyTotik1}.
Now apply (\ref{yut}) to this $\g$, to $z_0$ and to the poles $\pm \sqrt {a_i}$
with the associated orders $n_i$, but for the $2k$-th derivative.
We get a rational function $U_{2n}$, $n=n_0+\cdots+n_m$,
with poles at $\pm \sqrt{a_i}$ of order at most $n_i$ such that $\|U_{2n}\|_\g\le 1$
and
\[ |(U_{2n}(z))^{(2k)}(0)|\ge (1-\e)^{2k}(1-o(1))
\left(\sum_{i=0}^m n_i\left\{\frac{\partial g_{\C\setm \g}(0,\sqrt {a_i})}{\partial {\bf n}_+}
+\frac{\partial g_{\C\setm \g}(0,-\sqrt {a_i})}{\partial {\bf n}_+}\right\}\right)^{2k},
\]
which yields, in view of (\ref{sdf}),
\be  |(U_{2n}(z))^{(2k)}(0)|\ge (1-\e)^{4k}(1-o(1))
\left(\sum_{i=0}^m n_i\left\{\frac{\partial g_{\C\setm \G^*}(0,\sqrt {a_i})}{\partial {\bf n}_+}
+\frac{\partial g_{\C\setm \G^*}(0,-\sqrt {a_i})}{\partial {\bf n}_+}\right\}\right)^{2k}.
\label{yut1}\ee
Note also that, by the maximum principle,
we have $\|U_{2n}\|_{\G^*}\le\|U_{2n}\|_{\g}\le 1$
because all the poles of $U_{2n}$ lie outside $\g$.

By the symmetry of $\G^*$ and of the system $\{\pm\sqrt{a_i}\}$  onto the origin then $U_n(-z)$ also has
this property, furthermore  $(U_{2n}(-z))^{(2k)}(0)=(U_{2n}(z))^{(2k)}(0)$,
so if we set ${\cal R}_{2n}(z)=\frac12(U_{2n}(z)+U_{2n}(-z)),$ then ${\cal R}_{2n}$ is an even
rational function for which (\ref{yut1}) is true if we replace in it
$U_{2n}(z)$ by ${\cal R}_{2n}(z)$. But then there is a rational function
$R_n$ such that ${\cal R}_{2n}(z)=R_n(z^2)$, and for this $R_n$ we have
that
$\|R_n\|_\G=\|{\cal R}_{2n}\|_{\G^*}\le 1$, and (see (\ref{yut1}))
\[(R_n(z^2))^{(2k)}(0)
\ge (1-\e)^{4k}(1-o(1))2^{2k}\left(\sum_{i=0}^mn_i\O_{a_i}(A)\right)^{2k}\]
where we used the equality of the two quantities in
(\ref{b36}) and (\ref{fffh}).
Note also that this $R_n$ has poles at $a_0,\ldots,a_m$ of orders at most
$n_0,\ldots,n_m$.
Now the argument used in the proof of Theorem \ref{thMarkovk}
via the Fa\`{a} di Bruno's formula shows that the preceding inequality is the same
as
\[|R_n^{(k)}(0)|
\ge (1-\e)^{4k}(1-o(1))\frac{2^{k}}{(2k-1)!!}\left(\sum_{i=0}^mn_i\O_{a_i}(A)\right)^{2k}.\]

A similar construction can be done for the other endpoint $B$ of $\G$,
and by taking the larger of the two expressions in these lower estimates
we finally conclude the last statement in Theorem  \ref{thMarkovk}
regarding the sharpness of (\ref{mainMarkovk}).\endproof
\bigskip

\noindent {\bf Acknowledgement.}
This research was partially done
while the first author had a postdoctoral position at the Bolyai
Institute, University of Szeged, supported by the European Research Council Advanced
Grant No. 267055.
He was also partially supported by the Russian Science Foundation under grant 14-11-00022
(sharpness
of Bernstein type inequalities).

The second author was supported
by the Bolyai Scholarship of the Hungarian
Academy of Sciences.

The third author was supported by NSF grant DMS 1564541.

The third (senior) author also wants to acknowledge that this research has been started by the
other two authors. In particular, they had Theorems \ref{thcurve} and \ref{tharc}
for analytic
curves and arcs (Theorem \ref{thcurve} by a somewhat different proof), and the third author has joined the
research at that point.

\bigskip

Sergei Kalmykov

Department of Mathematics, Shanghai Jiao Tong University

800 Dongchuan RD, Shanghai, 200240, China

and

Far Eastern Federal University

8 Sukhanova Street, Vladivostok, 690950, Russia

{\it sergeykalmykov@inbox.ru}

\medskip{}

B\'ela Nagy

MTA-SZTE Analysis and Stochastics Research Group

Bolyai Institute, University of Szeged

Szeged, Aradi v. tere 1, 6720, Hungary

{\it nbela@math.u-szeged.hu}

\medskip

Vilmos Totik

MTA-SZTE Analysis and Stochastics Research Group

Bolyai Institute, University of Szeged

Szeged, Aradi v. tere 1, 6720, Hungary

\smallskip
and
\smallskip

Department of Mathematics and Statistics, University of South Florida

4202 E. Fowler Ave, CMC342, Tampa, FL 33620-5700, USA
\medskip

{\it totik@mail.usf.edu}
\end{document}